\documentclass{article}
\usepackage{epsfig}
\usepackage{bbm}
\usepackage{graphicx}
\usepackage{amsmath,amssymb}

\usepackage{booktabs}
\usepackage{makeidx}
\newcommand{\comment}[1]{}
\newcommand{\mathi}{\textrm{i}}
\newcommand{\EEA}{\end{eqnarray}}
\newcommand{\BEA}{\begin{eqnarray}}

\newcommand{\lh}{\hat{\lambda}}
\newcommand{\wb}{\omega_b}
\newcommand{\lb}{\lambda_b}

\newcommand{\gh}{\hat{\gamma}}
\newcommand{\btilde}{\tilde{b}}
\newcommand{\gtilde}{\tilde{\gamma}}
\newcommand{\sg}{\sigma_\gamma}
\newcommand{\dg}{d_\gamma}
\newcommand{\Sb}{\sigma_b}
\newcommand{\gb}{\gamma_b}
\newcommand{\sx}{\sigma_x}

\newcommand{\peq}{p_{eq}}
\newenvironment{example}[1][Example]{\begin{trivlist}
\item[\hskip \labelsep {\bfseries #1}]}{\end{trivlist}}
\usepackage[bookmarks,colorlinks,breaklinks]{hyperref}
\hypersetup{
    colorlinks,%
    citecolor=black,%
    filecolor=black,%
    linkcolor=black,%
    urlcolor=black
}

\makeindex

\begin{document}
\title{\bf Model error in data assimilation\footnote{The note is prepared for a chapter in ``Nonlinear and Stochastic Climate Dynamics",~Eds: C.L.E.~Franzke~and T.J. O'Kane, Cambridge University Press.}}
\author{John Harlim\footnote{email: jharlim@psu.edu.} \\ Department of Mathematics and Department of Meteorology \\ The Pennsylvania State University}
\date{\today}
\maketitle

Data assimilation\index{data assimilation} (or Bayesian filtering)\index{Bayesian filtering} is a statistical method to find the conditional distribution of the hidden variables of interest given noisy observations from nature. In application, the hidden variables of interest can be the state variables that are directly or indirectly observed or can even be some unobserved parameters in the models. In practice, data assimilation is typically realized by numerical schemes that produce conditional statistics of the state variables of interests, accounting for the information from the observations, rather than the corresponding conditional distribution; this gives a reasonable justification why we called it a ``statistical method". When observations are available at discrete times, Bayesian filtering is an iterative predictor-corrector scheme that adjusts the  prior forecast (background) statistical estimates from a predictor (or dynamical model) to be more consistent with the current observations. This correction step is referred to as analysis in the atmospheric and ocean science (AOS) community. Subsequently, the posterior (corrected or analysis) statistical estimates are fed into the model as initial conditions for future time prior statistical estimates. 

While the typical problems of interest are nonlinear, non-Gaussian, and high-dimensional, many practical data assimilation schemes that are currently used rely on Gaussian and/or linear assumptions. In particular, most practical data assimilation schemes are some type of approximation of the celebrated Kalman filter \cite{kalman:61}\index{Kalman filter}, which is the optimal solution (in the least squares sense) of the Bayesian filtering problem under linear and Gaussian assumptions. Essentially, all of these approximations were introduced to reduce the computational cost and to improve the statistical predictions. For example, in the AOS data assimilation community, two important schemes are: {\bf (i)} the ensemble Kalman filtering methods \cite{evensen:94,houtekamer:98,anderson:01,whitaker:02,anderson:03,bishop:01,evensen:03,hunt:07,szunyogh:05}\index{ensemble Kalman filter} which rely on empirical statistical estimates from ensemble forecasts; {\bf (ii)} variational-based methods \cite{sasaki:70,lorenc:86,courtier:94,lorenc:03} that rely on linear tangent and adjoint models. Operationally, most of the weather prediction centers, including the European Center for Medium-range Weather Forecasts (ECMWF), the UK Met Office, and the National Centers for Environmental Prediction (NCEP), are adopting hybrid approaches, taking advantages from both the ensemble and variational based methods \cite{ibbfhlr:10,bri:11,big:12,clb:13,wpkw:13}. While these approximations were introduced for practical consideration, theoretical understanding of the convergence of these methods in idealistic settings were established \cite{mcb:11,gth:13,km:15} for ensemble Kalman filter and for variational based methods \cite{lss:14,bllmss:13,blsz:13}. We should also mention that while these approximate methods can provide reasonable estimates of the first-order statistics, recent study \cite{ls:12} suggested that one should be cautious in interpreting their second order statistical estimates. 

In the numerical weather forecasting applications, these approximate filters are routinely used to assimilate observations collected from aircraft, radiosonde, satellite, and radar measurements to provide initial conditions for the atmospheric weather models, known as GCM's (General Circulation Models)\index{GCM (General Circulation Model)}. Depending on the model resolutions, the current GCM's have state variables of dimension $\mathcal{O}(10^9-10^{12})$. While model improvement is a significant enterprise that is continuously exercised, model error is unavoidable. This problem is attributed to incomplete understanding of the underlying physics and our lack of computational resources to resolve physical processes at various time and length scales or to model the interaction across scales. In the context of numerical weather prediction, the model is more accurate in the midlatitude atmospheric region since the dynamics can be approximated by quasi-geostrophic models that are well understood. In the Tropics, this approximation is not adequate since the Coriolis force vanishes at the equator and the dynamics is dominated by vertical heating/cooling in response to diabatic heating caused primarily by latent heat release. Despite some improvement in tropical weather forecasting \cite{bkjdrlrvb:08}, the forecast error for the zonal (east-west direction) wind component remains the largest in the Tropics (e.g., see Fig~1 in \cite{zitt:13}). 
The difficulty in predicting the Tropics is primarily caused by limited representation of the tropical convection and its multiscale organization in the contemporary convection parameterization \cite{moncrieff:07,bkjdrlrvb:08}. This is an example of  ``intrinsic information barrier" that prevents one from capturing the large-scale phenomena with a coarse model, as pointed out in \cite{mg:11b}.

Given these practical issues, an important challenge in data assimilation is to intelligently utilize the existing methods (either the ensemble, variational, and any hybrid based approaches) in the presence of model error. One difficulty is that model error can arise from any sources, such as imprecisedly specified model parameters, boundary conditions, unresolved processes, numerical approximation, etc. While the overall goal is to understand the implication of model error of any type on data assimilation, we emphasize on the effect of model error from unresolved scales. Our choice is partly because model error from unresolved scales is a subject of interest in applied mathematics under different names and a long list of literature exists in this subject. One goal of this chapter is to discuss this challenging issue in the simplest possible setting to clarify the problem and, subsequently, review some of the existing approaches to mitigate this issue. We will classify the existing approaches into two groups: those that estimate lower-order model error statistics directly are grouped as the \emph{statistical methods}; those that implicitly estimate the model error statistics with stochastic models beyond unbiased Gaussian white noises are grouped as \emph{stochastic parameterization methods}.\index{stochastic parameterization} We will discuss the pro and cons of all of these methods in the most transparent manner with simple examples. Subsequently, we use illuminating examples to understand the theory behind filtering with model error from unresolved scales that was recently established in \cite{bh:14}. This theory will justify why stochastic parameterization, as one of the main theme of this book, is an adequate tool for mitigating model error in data assimilation. Subsequently, we also discuss the main challenges in implementing stochastic parameterization in general. We will briefly discuss a recently proposed semiparametric framework as an alternative approach to mitigate these challenges \cite{bh:15qjrms}. This data-driven framework implements the nonparametric diffusion forecasting models \cite{bgh:15,bh:15physd} to represent dynamically evolving parameters in the existing data assimilation framework. We close this chapter with a short summary.

\section{Model error: a prior distribution formulation}\label{sec1}\index{prior distribution formulation}

Classical approaches to mitigate model error in data assimilation are motivated by analyzing the moments of the difference between the prior forecast estimate and the truth. To clarify this statement, we define model error\index{model error} as the difference between the truth $x(t)$ and the estimate $\tilde{x}(t)$ from imperfect model,
\begin{align}
e(t) \equiv x(t)-\tilde{x}(t).\label{modelerror}
\end{align}
In this context, we assume that the error, $e(t)$ is a random process, where the randomness can be due to the uncertainties in the initial error and/or the chaotic nature of the truth and the stochastic nature of the estimates. Even if the model is deterministic, the stochastic nature of the estimates is obvious when the estimates are outcomes of assimilating noisy observations. 

We assume that the random process in \eqref{modelerror} has mean $\bar{e}(t)=\mathbb{E}[e(t)]$ and covariance $Q(t)=\mathbb{E}[(e(t)-\bar{e}(t))(e(t)-\bar{e}(t))^\top]$. Similarly, we define $\bar{x}(t)\equiv \mathbb{E}[x(t)]$  and $P(t)= \mathbb{E}[(x(t)-\bar{x}(t))(x(t)-\bar{x}(t))^\top]$ as the forecast mean and covariance estimates from the perfect model (i.e., the true mean and covariance statistics), respectively. 
We also define $\bar{\tilde{x}}(t)\equiv\mathbb{E}[\tilde{x}(t)]$ and $\tilde{P}(t)=\mathbb{E}[(\tilde{x}(t)-\bar{\tilde{x}}(t))(\tilde{x}(t)-\bar{\tilde{x}}(t))^\top]$, as the prior mean and covariance estimates from the imperfect model. One can show that the mean model error,
\begin{align}
\bar{e}(t) = \bar{x}(t)-\bar{\tilde{x}}(t),\label{modelerrormean}
\end{align}
is essentially the ``bias forecast error"\index{bias}, defined in \cite{dds:98}. Taking the expectation square of the difference between  $x(t)$ in \eqref{modelerror} and $\bar{x}(t)$ in \eqref{modelerrormean}, one can deduce the forecast error covariance,\index{forecast error covariance}
\begin{align}
P(t) &= \tilde{P}(t) + \Big(Q_{\tilde xe}(t) +Q_{e\tilde x}(t) +Q(t)\Big),\label{modelerrorcov}
\end{align}
where $Q_{\tilde xe}(t)= \mathbb{E}[(\tilde x(t)-\bar{\tilde x}(t))(e(t)-\bar{e}(t))^\top]$ and $Q_{e\tilde x}(t)=Q_{\tilde xe}^\top(t)$ denote the cross covariances between the forecast state, $\tilde{x}(t)$, from the imperfect model and the model error estimator $e(t)$. Equations~\eqref{modelerrormean} and \eqref{modelerrorcov} suggest that in the presence of model error the first two-order statistics of the truth can only be recovered when the bias, $\bar{e}(t)$, are added to the prior mean estimates $\bar{\tilde{x}}(t)$ and the prior error covariances, $\tilde{P}(t)$, are appropriately adjusted by covariance correction factors $Q_{\tilde xe}(t)+ Q_{e\tilde x}(t)+ Q(t)$. One can obviously repeat this formalism on higher-order moments of interest but they are not important in our discussion here. 

The most important fact that one should realize behind this implicit formalism is that while the formula looks deceptively simple, it does not provide any easy access to the model error statistics, even for the lower order statistics such as $\{\bar{e}, Q_{\tilde xe}, Q\}$ in \eqref{modelerrormean} and \eqref{modelerrorcov}. It is worthwhile to point out that even if we know the dynamics of $e(t)$, while the problem becomes simpler, its statistics may not be easily determined explicitly in practical situation. To see this, suppose the joint variables $(\tilde{x},e)$ solve a system of differential equations,
\BEA
\frac{d\tilde x}{dt} = f(\tilde x, e), \quad
\frac{de}{dt} = g(\tilde x, e),\nonumber
\EEA
where for simplicity $f$ and $g$ are assumed to be deterministic and known. Assume also that $(\tilde{x},e)$ can be characterized by a joint density function $p(\tilde x, e,t)$, which solves the corresponding Liouville equation \cite{lm:94}\index{Liouville equation}, $\partial_tp = -\nabla_{\tilde{x}}\cdot (fp)-\nabla_e\cdot (gp)$. Then the model error statistics \index{model error statistics} solve a system of differential equations for the moments of the Liouville equation. In general, however, these differential equations can be infinite-dimensional since the moments may interact with all of the higher order moments. To see this, consider the following simple example.
\begin{example}[Example 1:]
Let us assume that the dynamics of $e$ is independent of $\tilde{x}$ and our aim is to compute the mean model error, $\bar{e}$.
Consider a simple model error estimator $e(t)\in\mathbb{R}$ that satisfies, 
\BEA
\frac{de}{dt} = g(e) = ae +be^2,
\EEA
for some constants $a, b$. Assume that $e(t)$ can be characterized by a marginal density function $p(e,t)$ that decays to zero as $e\rightarrow\pm\infty$ and that it solves the Liouville equation:
\BEA
\frac{\partial p}{\partial t} = - \frac{\partial}{\partial e}[gp] = - \frac{\partial}{\partial e}[(ae+be^2)p].\label{liouville}
\EEA
The first moment can be computed by multiplying \eqref{liouville} with $e$ and taking an expectation (or integral with respect to $\mathbb{R}$) such that,
\BEA
\frac{d}{dt} \int_\mathbb{R} e p\,de &=& \int_\mathbb{R} e \frac{\partial p}{\partial e}\, de = - \int_\mathbb{R} e\frac{\partial}{\partial e}[(ae+be^2)p]\,de = \int_\mathbb{R} (ae+be^2)p\,de, \nonumber \\ 
&=& \int_{\mathbb{R}} (ae +2be\bar{e} -b\bar{e}^2 + b(e-\bar{e})^2) p\,de
\label{dmean}
\EEA
where we use integration by part and the standard completing square trick.
Since $\bar{e}=\mathbb{E}[e] = \int_\mathbb{R} e p\,de$ and $Q = \mathbb{E}[(e-\bar{e})^2] = \int_\mathbb{R} (e-\bar{e})^2 p\,de$, we can rewrite \eqref{dmean} as follows,
\BEA
\frac{d\bar{e}}{dt} = a\bar{e} + b\bar{e}^2 + bQ = g(\bar{e})+ bQ.\nonumber
\EEA
The differential equation for the variance can be deduced by multiplying \eqref{liouville} with $(e-\bar{e})^2$ and taking an integral with respect to $\mathbb{R}$. Repeating the same algebraic manipulation, we obtain,
\BEA
\frac{dQ}{dt} = 2g'(\bar{e})Q + 2bS. \nonumber
\EEA
where we define $S\equiv\mathbb{E}[(e-\bar{e})^3]$ as the third-order centered moment. Notice that in this very simple example, the mean model error $\bar{e}$ depends on $Q$, the model error covariance $Q$ depends on $S$, and one can check that the $S$ will depend on higher-order moments. Therefore, the system of differential equations for the dynamics of the statistics of $p$ is not closed. This issue reminisces the classical turbulent closure problem, where the expectation is replaced with the Reynold averaging\index{Reynold averaging}. While the example is only one-dimensional, the computational costs significantly increase for high-dimensional model error estimator $e(t)$, even if we apply some higher-order moment truncation. Essentially, if $e\in\mathbb{R}^n$, then $Q$ has $n^2$ components, $S$ has $n^3$ components, and so on. 
\end{example}
 
In real application, the problem is much more difficult since we have no access to either the truth $x(t)$ or the model for $e(t)$. In the next section, we discuss some of the existing methods for estimating the model error mean and covariance statistics.  

\section{Estimating model error statistics in data assimilation} \index{model error statistics}

Since most data assimilation methods that are used in the numerical weather forecasting application produce mean and covariance estimates (except for the variational based methods that are usually implemented with a fixed covariance matrix), then mitigating model error is highly associated to finding the model error mean, $\bar{e}$, and covariance statistics, $Q_{\tilde xe}$ and $Q$. In particular, given the analysis mean, $\bar{x}^a_{m-1}$, and covariance estimate, $P^a_{m-1}$, at a particular instance $t_{m-1}$, one uses the (imperfect) model to propagate these statistics to obtain $\bar{\tilde{x}}^b_{m}$ and $\tilde{P}^b_{m}$ at the next observation time, $t_m$. Subsequently, one uses the estimated model error statistics $\{\bar{e}(t_m), Q_{\tilde xe}(t_m), Q(t_m)\}$ to adjust the prior statistical estimates $\bar{\tilde x}^b_{m}$ and $\tilde{P}^b_{m}$ to be closer to the corresponding true prior statistics, $\bar{x}^b_m$ and $P^b_m$, respectively. To close the cycle, one applies his/her data assimilation method of choice to obtain the analysis mean, $\bar{x}^a_m$, and covariance estimate, $P^a_m$, accounting for observations at the current time. 

In this section, we discuss several practical methods that directly estimate the mean and covariance statistics of the model error estimator, $e(t)$. 

\subsection{Classical state-augmentation approach}\label{sec21}\index{state-augmentation method}
The simplest type of model error is the misspecification of constant parameters in the dynamical model. For this scenario, one can apply various statistical methods to estimate the unknown parameters. We refer to this type of model error as the simplest in the sense that the source of model error is known, that is, through a misspecification of constant parameters. By simplest here, we do not say that the parameters are easily estimated; this will depend on the identifiability of the parameters and the estimation schemes. When the parameters are time dependent, the complexity of the problems is significantly increased. 

A classical approach for estimating parameters is to apply Kalman filter based methods on an augmented vector of state-parameter \cite{friedland:69,friedland:82},  
\begin{equation}\label{stateaugmented}
\begin{aligned}
\dot{\tilde x} &=& f(\tilde x,\theta),\\
\dot{\theta} &=& g(\tilde x,\theta),
\end{aligned}
\end{equation}
where ``the dot" indicates time derivative and $g$ is typically chosen empirically. If the correct parameters are obtained, the model error vanishes and the true prior statistics are directly obtained. When the true parameters $\theta$ are constant and identifiable, one can just apply this strategy with persistence model, $g=0$.
The main issue with this approach is that if the parameters are not constant, then choosing the appropriate parametric model, $g$, can be difficult; classical choices often assumed $g$ to be independent of $\tilde x$. Later in section~\ref{sec5}, we will discuss a nonparametric approach for modeling $\theta$, assuming that it is independent of $\tilde x$ \cite{bgh:15}. 

While the state augmentation approach was designed as a parameter estimation method, it can also be implemented to directly estimate the mean model error or forecast bias \cite{dds:98}; that is, solving the augmented equations in \eqref{stateaugmented} for $\bar{x}$ and $\bar{e}$ (in place of $\tilde{x}$ and $\theta$, respectively). In practice, this procedure was implemented with an additional assumption for the model error covariance. The typical choice is to assume the model error covariance to be proportional to the forecast error covariance \cite{dds:98}, that is,
\BEA
Q(t_m)\approx \alpha \tilde{P}^b_m,\label{Qadhoc} 
\EEA
where $\tilde{P}^b_m$ is the prior covariance estimate from imperfect model and $\alpha$ is an empirically chosen scalar. Technically, this approach estimates the model error covariance $Q$ with a multiplicative covariance inflation of $\tilde{P}^b_m$ and ignores the cross covariances $Q_{\tilde{x}e}$ completely. Moreover, the parameter model $g$ is often chosen on an ad hoc basis, such as the persistence model\index{persistence model}, $g=0$, or the white noise processes, $g=\sigma \dot{W}_t$, \cite{friedland:69,friedland:82}, with an empirically chosen noise amplitude $\sigma$. We should note that other models for $g$ were also proposed in \cite{baek:06} with empirical choices of $Q$. 

\subsection{Estimating model error covariances}\label{sec22}

While estimation of the forecast bias term, $\bar{e}$, is important in mitigating model error, many approaches to account for model error are heavily associated with estimating matrix $Q$ in \eqref{modelerrorcov}. This traditional point of view is often based on assuming that model error can be treated as unbiased Gaussian white noise processes \cite{tb:88,vt:02,kalnay:03},\index{unbiased Gaussian white noise processes}
\BEA
e(t) = x(t) - \tilde{x}(t) \approx \eta(t), \quad \eta(t)\sim\mathcal{N}(0,Q(t)),\label{whitenoise}
\EEA
where in some applications, further stationarity assumption may take place by setting $Q$ to be time-independent. While this tacit assumption can be useful for some problems, it may not always produce satisfactory estimates since model error also introduces forecast bias $\bar{e}$, cross covariance statistics, $Q_{\tilde xe}$ and other higher-order statistics that are needed for accurate statistical estimation.

Just to name a few examples, in the weak constrained 4D-VAR implementation \cite{tremolet:07}, unbiased model error is assumed, that is, $\bar{e}=0$, in addition to an empirically chosen model error covariance estimator. Typical choice sets $P(t)=B$ in \eqref{modelerrorcov}, for a fixed covariance matrix $B$. In this case, $Q_{\tilde xe}, Q$ and $\tilde{P}$ are implicitly accounted and the accuracy of estimates relies on the methodology for choosing the $B$ matrix. In the ensemble Kalman filter community,\index{ensemble Kalman filter} such practice (setting $\bar{e}=0$ and modeling error covariance with \eqref{Qadhoc}) is known as ``multiplicative covariance inflation"\index{multiplicative covariance inflation}; this practical approach was introduced to mitigate  covariance underestimation 
due to unresolved scales model error \cite{hw:05,wh:12} or when small ensemble size is used \cite{aa:99}. An alternative approach known as  ``additive covariance inflation"\index{additive covariance inflation} was also used to account for inhomogeneity of the underestimated covariance matrix \cite{ott:04,scyangetal:06,kalnayetal:07,whwst:08}. In practice, one prefers the multiplicative covariance inflation rather than the additive covariance inflation since it is difficult to specify the appropriate ansatz for the additive inflation matrix with appropriate scaling when the system variables have different quantifying units (personal communication with J.L. Anderson). There is also a relaxation-to-prior method \cite{zss:04} that was found to be useful in various applications; this method adjusts the analysis error covariance to be closer to its prior error covariance estimate with an empirical chosen adjustment coefficient. Note that this approach implicitly approximates $Q_{\tilde xe}$ with the empirical cross-covariances between the prior and posterior errors, in addition to $Q$. A more systematic Bayesian approach that alleviates the covariance undersampling in the ensemble Kalman filter context was studied by \cite{bouquet:11}.  

Alternatively, adaptive methods to estimate covariance statistics have been proposed since early 70's \cite{mehra:70,mehra:72,belanger:74}. These methods were designed to estimate covariance matrices, $Q, R,$ and $C$ of the following discrete-time linear stochastic filtering problem,
\BEA
x_{m+1} &=& Fx_m + e_m, \quad e_m \sim\mathcal{N}(0,Q),\\
v_m &=& Hx_m + \sigma_m, \quad \sigma_m \sim\mathcal{N}(0,R),\label{noisyobs}
\EEA
where $F$ and $H$ are linear dynamical and observation operators, respectively. The formulation also allows one to estimate $\mathbb{E}(e_m\sigma_m^\top)=C$, the correlation between the system and observation error noises. In this filtering problem, the truth is stochastic; it is inherently driven by unbiased white noise processes. Therefore, the covariance matrix $Q$ is not associated with any model error term. 

The main idea of these adaptive covariance estimation methods\index{adaptive covariance estimation methods} \cite{mehra:70,mehra:72,belanger:74} is to apply Bayes' theorem\index{Bayes' theorem} to obtain a posterior distribution of the augmented state and parameters at each time step $t_m$ when observations become available,  
\begin{align}
p(x_m,\theta_m| v_m)&\propto p(x_m,\theta_m) p(v_m| x_m,\theta_m),\label{bayes}
\end{align}
where $p(x_m,\theta_m)$ denotes the prior distribution of the augmented state and parameters at time $t_m$ and $p(v_m| x_m,\theta_m)$ denotes the likelihood function of the augmented variables, corresponds to the observation model in \eqref{noisyobs}. In this formalism, $\theta_m =(Q(t_m),R(t_m),C(t_m))$. The parameterization method can be formally described as follows: Since $p(x_m,\theta_m) = p(\theta_m)p(x_m|\theta_m)$ by definition of the conditional distribution, we can rewrite \eqref{bayes} as follows:
\begin{align}
p(x_m,\theta_m| v_m)&\propto  p(\theta_m) p(x_m|\theta_m) p(v_m| x_m,\theta_m),\label{bayes2}\\
&\propto p(\theta_m) p(x_m|\theta_m,v_m).\label{bayes3}
\end{align}
Here, the first step in the filtering algorithm is to estimate $p(x_m|\theta_m,v_m)$ by applying Bayes' theorem to the last two components of \eqref{bayes2}. Subsequently, we implement the Bayes' theorem one more time in \eqref{bayes3} to obtain the posterior distribution of the augmented variables $(x_m,\theta_m)$.

Numerically, the two-step Bayes' update in \eqref{bayes2}-\eqref{bayes3} can be approximated with the Kalman filter \cite{mehra:70,mehra:72,belanger:74} or an extended Kalman filter for nonlinear systems in \cite{dcdg:85}. We should mention that the paper \cite{dcdg:85} also provides an efficient implementation for the method in \cite{belanger:74}. Recent extension of these methods using ensemble Kalman filters were shown in \cite{bs:13,hmm:14,zh:15}.\index{ensemble Kalman filter} Indeed, numerical comparisons between the three methods in \cite{bs:13,hmm:14,zh:15} were shown on various examples in \cite{zh:15}. Practically, the first step is to apply a primary filter (either KF, EKF, or EnKF) to update the statistics for $x_m$, and subsequently, a secondary filter is used to update $\theta_m$. In the secondary update, the prior model for $p(\theta_m)$ is typically empirically chosen, e.g., the persistence model, $\theta_{m+1}= \theta_m$. To avoid unobservability of the parameters due to sparse observations with dimension less than the number of parameters, $\theta_m$, one includes information from past observations up to lag $L>1$ (see \cite{belanger:74,hmm:14,zh:15} for methods that can use $L\geq 1$).  

Practically, the primary filter produces Gaussian statistics of $p(x_m|\theta_m,v_m,\ldots,v_{m-L+1})$ since Kalman-based filters are used in estimating $x_m$. However, the dependence of $p$ on $\theta_m$ can be described non-uniquely \cite{mehra:70,mehra:72,belanger:74,bs:13,hmm:14,zh:15}. For example, Belanger's formulation \cite{belanger:74} defines $p(x_m|\theta_m,v_m,\ldots,v_{m-L+1})$, as a likelihood function of $\theta_m$, through the following observation model,
\begin{align}
\sigma_{m,\ell} =  \mathcal{F}_{m,\ell}\theta_m + \eta_{m,\ell},\quad \eta_{m,\ell}\sim\mathcal{N}(0,W_{m,\ell}), \quad \ell =m,\ldots, m-L+1,\label{pobsmodel}
\end{align}
for any lags $L\geq 1$.  
In \eqref{pobsmodel}, components of $\sigma_{m,\ell} = \{d_m d_{m-\ell}^\top\}$ are the product of the forecast error estimates in the observation space,
\begin{align}
d_m = v_m - H\bar{x}^b_m,
\end{align}
where $\bar{x}^b_m$ denotes the mean prior estimate obtained from the primary filter. In \eqref{pobsmodel}, the observation operator $\mathcal{F}_{m,\ell}$ and the noise covariance matrix $W_{m,\ell}$ are functions of $\bar{x}^b_{m-\ell}$ and they will be constructed recursively. We should also note that $W_{m,\ell}$ is typically approximated under Gaussian assumption (see \cite{belanger:74,hmm:14} for the detail formula of $\mathcal{F}_{m,\ell}$ and $W_{m,\ell}$). With the pseudo-observation model in \eqref{pobsmodel}, a secondary Kalman filter is implemented $L$-times to sequentially update the posterior mean and covariance estimate of $\theta_m$, accounting for the pseudo-observations $\{\sigma_{m,\ell}\}_{\ell=1\ldots, L}$ one at a time. We should note that there are other methods to approximate the secondary Bayes' update in \eqref{bayes3} that use different observation model  and do not use Kalman update, see e.g., \cite{bs:13,zh:15}.  

While these adaptive covariance estimation methods were not designed to estimate $Q$ associated with model error, it can be used to estimate the model error covariance $Q$, assuming that the model error estimator is an unbiased white noise process as in \eqref{whitenoise}. In this particular application, the covariance estimation method essentially acts like an adaptive covariance inflation\index{adaptive covariance inflation} method of an additive type. We should mention that while adaptive covariance inflation methods have been proposed \cite{anderson:07,lkm:09,miyoshi:11}, they are all multiplicative type; they adaptively estimate the multiplicative factor $\alpha$ in \eqref{Qadhoc}. In the following, we will demonstrate a numerical example showing application of the adaptive covariance estimation method as an additive adaptive covariance inflation method for ETKF\index{ETKF (Ensemble Transform Kalman Filter)} \cite{hunt:07} to mitigate model error. 

\begin{example}[Example 2:]
Consider a data assimilation experiment with the Lorenz-96 model \cite{lorenz:96},\index{Lorenz-96 model}
\BEA
\frac{dx_j}{dt} = (x_{j+1}-x_{j-2})x_{j-1} - \theta x_j + 8,\label{L96}
\EEA
with an additional parameter $\theta$. In \eqref{L96}, index $j=1,\ldots,40$ represents the spatial grid point with periodic boundary condition. Let the truth be solutions of \eqref{L96} with $\theta=1$ at every time interval $\Delta t=0.05$ (which corresponds to the standard Lorenz's 6-hour time unit). Suppose that model error is committed from specifying $\theta=1.2$ which is different than the true value of $\theta=1$ but this misspecification is unknown. Suppose noisy observations of $x_j$ at every grid point are collected; these observations are corrupted with i.i.d.~Gaussian noises with mean zero and \emph{unknown} error variance, resulting to a 40-dimensional identity error covariance matrix $R=\mathcal{I}_{40}$ without cross correlation $C=0$.

We will now employ the ETKF method (based on the formulation in \cite{hunt:07}) with the additive adaptive covariance inflation method discussed above. The detail of the algorithm is in the Appendix of \cite{hmm:14}. In this implementation, we implicitly assume that the model error, $e(t)$, is modeled as unbiased white noise Gaussian processes as in \eqref{whitenoise} with a covariance structure,
\BEA
Q = \begin{pmatrix} q_1 & q_2 & 0 & 0 & \ldots & 0 & q_2 \\
q_2  & q_1 & q_2 & 0 & \ldots & 0 & 0 \\ 
0 & q_2  & q_1 & q_2 & \ldots & 0 & 0 \\
\vdots & & \ddots & \ddots & \ddots & & \vdots \\
q_2 & 0 & 0 & \ldots & 0 & q_2 & q_1 
\end{pmatrix},
\EEA
so we will estimate only two parameters $\{q_1,q_2\}$ for $Q$. Of course, this choice is adhoc, and it is partially motivated by the isotropic characteristic of the Lorenz-96 model in \eqref{L96}. Different choices of $Q$ were used in \cite{bs:13,zh:15}. In addition to these two parameters, we also estimate the variance of $R=r\mathcal{I}_{40}$, which true value is $r=1$.

In Figure~\ref{figL96}, we compare the results with two ensemble sizes, 10 and 20. Both experiments use ensemble sizes that are considerably smaller than the model state space $40$. Obviously, the larger ensemble size experiment produces much better results (the absolute error for the $10$-th component of $x$ is smaller than that with smaller ensemble size). Notably, the absolute error for the observation error covariance parameter, $r$, is closer to zero and smaller covariance inflation parameters $q_1,q_2$ are obtained in the experiment with ensemble size of 20. When ensemble size is smaller, 10, sampling error becomes more severe so the adaptive filter weights more to the observations by reducing the estimate for $r$ below its true value (with larger absolute error in $r$). The net effect of this reduced estimate in $r$ is similar to implementing a multiplicative covariance inflation. Simultaneously, the adaptive filter implicitly applies larger additive covariance inflation with much larger $q_1, q_2$. We suspect that one can improve this result with appropriate choices of $Q$; here, our goal is only to demonstrate that an adaptive covariance estimation method can be used as an adaptive covariance inflation method in the presence of model error. 

\begin{figure*}
\centering
\includegraphics[width=0.7\textwidth]{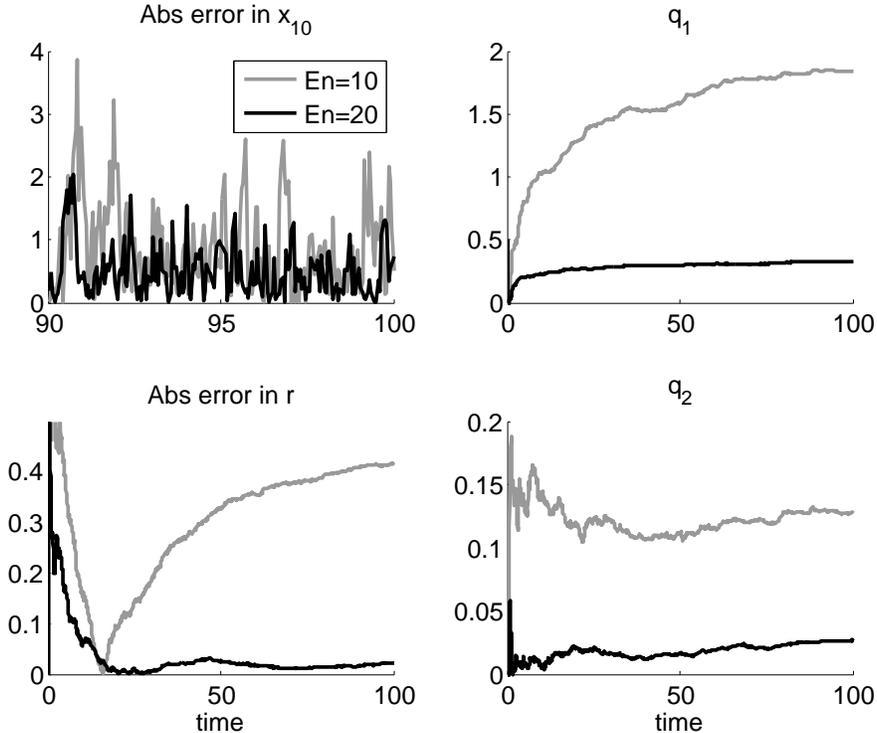}
\caption{Absolute errors of the posterior mean estimates (grey for ensemble size 10 and black for ensemble size 20) for $x_{10}$ [left upper panel] as functions of time. Absolute errors of $r$ as functions of time [left lower panel]; as a reference, the true $r$ is one. The remaining two panels depict the trajectories of the parameters $q_1, q_2$.}
\label{figL96} 
\end{figure*}
  
\end{example}

\subsection{Simultaneous estimations of bias and model error covariances}

The numerical approaches discussed above estimate only one of the model error statistics, either the mean or covariance, imposing various assumptions on the other statistics that are not estimated. Furthermore, these methods ignore estimating the cross-covariances, $Q_{\tilde xe}$, between the prior forecast, $\tilde{x}(t)$ and the model error, $e(t)$.  

There are few adhoc methods that simultaneously estimate for, both, the forecast bias (or mean model error) as well as the model error covariance. For example, the multiphysics (or multimodel) ensemble approach\index{multimodel ensemble approach} was proposed for simulating surrogate statistics for the model error such as the forecast bias (see the review article \cite{mz:11}). Recent mathematical justification of such approaches was given through an information theoretic framework\index{information theoretic} \cite{bm:15}. Another method, proposed by \cite{lkmd:09}, simultaneously estimates a certain parametric form of mean model error estimator and applies an empirical choice of additive covariance inflation. 

Alternatively, there is also the ``deterministic approach" introduced by \cite{cvn:08,cv:10,cv:11}. Their approach is motivated by the deterministic formulation\index{deterministic formulation} that was introduced in \cite{nicolis:04} which approximates the model error dynamics with a short-time Taylor's expansion about the initial conditions. That is, suppose that $\dot{x}= f(x)$ denotes the true dynamics and $\dot{\tilde x}=\tilde f(\tilde x)$ denotes the imperfect model, then they approximate the model error dynamics as follows,
 \BEA
e(t) = \int_0^t (\dot{x}-\dot{\tilde x})d\tau= \int_0^t \Big(f(x(\tau)) - \tilde{f}(\tilde x(\tau))\Big) d\tau \approx \Big(f(x(0)) - \tilde{f}(\tilde x(0))\Big) t,
\EEA
employing a first-order Taylor's truncation. Subsequently, the model error mean and covariance statistics, including the cross covariance statistics $Q_{\tilde{x}e}$, are estimated by taking an empirical ensemble average locally about the initial conditions. In their implementation, they apply further approximations to obtain these statistics; either with the reanalysis data \cite{cvn:08} or some linear tangent approximation of $\tilde{f}$ \cite{cv:11}. Improvement for short-time assimilation interval were shown \cite{cv:11} and the method starts to diverge for longer time assimilation interval, by design.    

In the next section, we discuss methods that implicitly account for all the model error statistics through stochastic modeling for $e$ beyond the traditional white noise approximation in \eqref{whitenoise}. In particular, we discuss treatment for model error from unresolved scales. One shall see that while these approaches do not directly estimate the statistics of the model error and some of them utilizes the algorithms that were described in Sections~\ref{sec21} and \ref{sec22} to estimate parameters in the stochastic model for $e$.

\section{Stochastic parameterization of model error from unresolved scales}

Model error from unresolved scales has been an important mathematical subject under different names; e.g., model reduction, averaging, homogenization, Mori-Zwanzig formalism\index{Mori-Zwanzig formalism}, just to name a few. Loosely speaking, the mathematical perspective is to obtain an approximate model that involves only the coarse-grained variables of interest. For practical reasons, one is interested to use the reduced models\index{reduced models} in place of the full dynamics and therefore model error is committed by not resolving all of the variables beyond these coarse-grained variables of interest. 

An elegant way to formulate this problem is through the Mori-Zwanzig formalism \cite{zwanzig:61,mori:65,zwanzig:73}, which is just a way of rewriting a system of differential equations without assuming any scale separation. Suppose the full dynamics are governed by a system of ODE's,
\BEA
\frac{dx}{dt} &=& f(x,y), \label{slow}\\
\frac{dy}{dt} &=& g(x,y),\label{fast}
\EEA 
with initial conditions, $x(0)=x_0, y(0)=y_0$, and, for simplicity, we assume that $x$ is the coarse grained variables of interest and $y$ denotes the unresolved variables. Following the notation in \cite{gks:04}, one can rewrite the dynamics of $x$ in a closed form known as the \emph{Generalized Langevin Equation} (GLE),\index{GLE (Generalized Langevin Equation)}
\BEA
\frac{dx}{dt} &=& \bar{f}(x) + \int_0^t K(x(t-s),s)\,ds + n(x_0,y_0,t).\label{GLE}
\EEA
Here the first term is a Markovian\index{Markovian} dynamics, resulting from a projection of $f$ onto a function of $x$ only. This implies that the GLE in~\eqref{GLE} is non-unique, it depends on the choice of the projection operator (see \cite{mori:65,chk:00} for various different choices of projection operators). Particularly relevant to our case is to define the projection operator to be a conditional expectation given $x$ such that, $\bar{f}(x) = \mathbb{E}[f|x]$, as suggested in \cite{chk:00,gks:04}. The last two terms in \eqref{GLE} are consecutively known as the ``memory" term\index{memory term} (defined with a memory kernel $K$) to represent feedback from the unresolved scales and the orthogonal dynamics\index{orthogonal dynamics} that is treated as ``noise" that depends on the randomness of the initial condition, $y_0$. See e.g., \cite{gks:04,chk:00} for detail derivation or more general formulation of the GLE. While the formula in \eqref{GLE} is exact, realizing the last two terms in \eqref{GLE} is not easier than solving the full problem in \eqref{slow}-\eqref{fast}. Hence, model error is unavoidable when some approximation is applied to estimate the last two terms in \eqref{GLE}. 

Many methods were introduced to estimate the memory and noise terms in \eqref{GLE}, see e.g., \cite{ch:13} and the references therein. On the other hand, there are many methods that were not directly designed to estimate the right hand terms in \eqref{GLE} but they implicitly simulate these terms. Just to name few examples of such methods that can be valuable for data assimilation applications: (i) superparameterization\index{superparameterization} \cite{mg:13}; (ii) the reduced order modified quasilinear Gaussian algorithm \cite{sm:13}; (iii) the physics-constrained multilevel nonlinear regression model \index{physics-constrained multilevel nonlinear regression model}\cite{mh:13,hmm:14}; (iv) Markov chain type modeling \cite{cev:08,kbm:10}\index{Markov chain models}; (v) Heterogeneous Multiscale Methods-based reduced models \index{HMM (Heterogeneous Multiscale Methods)} \cite{eelrv:07}. (vi) Classical turbulent closure methods such as the Direct-Interaction Approximation (DIA)\index{DIA (Direct-Interaction Approximation)} for parameterizing subgrid scale processes in isotropic turbulence \cite{kraichnan:59} and its derivatives, Quasi-diagonal DIA, cumulant update DIA, and the regularized cumulant update DIA \cite{of:04} for modeling non-Markovian memory in inhomogeneous turbulence over topography. As of the authors knowledge, among all methods mentioned above, two of them that have made direct impact are the approach in (i) and (iv) which have been used to model cloud processes in GCM's \cite{grabowski:04,kscmt:11}. Some versions of (i) were also used for simulating combustion problems \cite{kerstein:88,kerstein:99}. The QDIA method in (vi) has also been proposed as a closure model for data assimilation \cite{of:10}.

While most of these approaches are derived from the first principle, assuming that the true dynamics are governed by a certain function (typically with simple prototype models of geophysical fluid dynamics), in practice, the modelers don't know what is the true dynamics. What available is a set of complex equations that is believed to be a reasonable approximation of the dynamics of $x$, such as the General Circulation Models (GCM's) in the weather and climate prediction community. In such a scenario, extension of these methods can be quite challenging. Conceptually, one way to formalize this issue is to assume that the available model, $\dot{\tilde{x}} = \tilde{f}(\tilde x)=\bar{f}(\tilde x)$, is an approximation of the unknown dynamics in \eqref{slow}-\eqref{fast}. Mathematically, this assumes that there exists a projection operator that maps $(x,y)$ to $x$, corresponding to the dynamical operator $\bar{f}$ given by the modelers. Of course, the projection operator is unknown in general and this assumption is only to give some intuition. A similar perspective was also described in \cite{ch:13}. Therefore, the model error estimator, $e(t) = x(t)-\tilde{x}(t)$, is non-Markovian,
\BEA
\frac{de}{dt} = \frac{dx}{dt}-\frac{d\tilde{x}}{dt} = \bar{f}(\tilde x + e) - \bar{f}(\tilde x) + \int_0^t K(\tilde x(t-s)+ e(t-s),s)\,ds + n(x_0,y_0,t),\label{edyn}
\EEA
expressed in terms of $\tilde x$ and $e$. If initial errors are zero, then the model error is intrinsically generated through the memory and noise terms in \eqref{edyn}. In this context, \emph{stochastic parameterizations}\index{stochastic parameterization} can be interpreted as methods to approximate this integro-differential equation. Three questions naturally arise: 
\begin{enumerate}
\item Which model should we use to approximate \eqref{edyn}? 
\item If the class of models to be used are in parametric form and Markovian, how do we estimate the parameters in these parametric models? 
\item Also, how do we ensure the stability of these parametric models?
\end{enumerate}
For the third question, it was shown that if the proposed parametric form is not carefully chosen, then the resulting model can blow up in finite time and gives no prediction skill \cite{my:12}. For a class of diffusion processes, a physics constrained parametric model has been proposed to overcome this issue \cite{mh:13,hmm:14}. While this strategy gives guidelines for choosing parametric models, it is indeed not easy to find one that produces accurate and consistent equilibrium statistical estimates as shown in the example in \cite{hmm:14} and hence we are back to question 1 above. We should note that the first two questions are wide open problems in general and are the same questions that have been posted in turbulent closure problems \cite{kraichnan:59}. Next, we review some numerical approaches below to gain some intuition.

\begin{example}[Example 3]\footnote{This example is taken from Section 4 of \cite{bh:14}}
Consider filtering the two-layer Lorenz-96 model \cite{lorenz:96}, whose governing equations are a system of $N(J+1)$-dimensional ODEs given by,\index{Two-layer Lorenz-96 model}
\begin{equation}\label{lor96}
\begin{aligned}\frac{dx_i}{dt} &= x_{i-1}(x_{i+1}-x_{i-2}) - x_i + F + h_x\sum_{j=(i-1)J+1}^{iJ} y_{j}, \\
 \frac{dy_{j}}{dt} &= \frac{1}{\epsilon}\big(a y_{j+1}(y_{j-1}-y_{j+2}) - y_{j} + h_y x_{\textup{ceil}(i/J)} \big),
\end{aligned}
\end{equation}
where $\vec{x}=(x_i)$ and $\vec{y}=(y_j)$ are vectors in $\mathbb{R}^{N}$ and $\mathbb{R}^{NJ}$ respectively and the subscript $i$ is taken modulo $N$ and $j$ is taken modulo $NJ$. In the example here, we set $N=8, J=32, \epsilon=.25, F=20, a=10, h_x = -0.4, h_y=0.1$. In this regime the time scale separation is small. To generate the observations, we integrate this model using the Runge-Kutta method (RK4) with a time step $\delta t=0.001$ and take noisy observations $\vec{v}_m\in\mathbb{R}^M$ at discrete times $t_m$ with various time intervals $\Delta t=t_{m+1}-t_m$ given by,
\begin{align}
\vec{v}_m = h(\vec{x}(t_m)) + \eta_m, \quad \eta_m\sim\mathcal{N}(0,R),\label{L96obs}
\end{align}
where $R=0.1\mathcal{I}_M$. In our experiment below, we set $M=4$ by taking observations at every other grid point. Suppose the reduced (or the available) model is the single layer Lorenz-96 model:
\BEA
\label{lor96_2}\frac{d\tilde x_i}{dt} = \tilde x_{i-1}(\tilde x_{i+1}-\tilde x_{i-2}) -\tilde x_i + F,
\EEA
such that the right-hand-term is $\bar{f}$ in notation \eqref{edyn} and let's try to address the two questions above. While systematic derivations to deduce the appropriate model error estimator for this simple model are available \cite{mtv:01,mfc:09}, they may be difficult to carry when the models are complicated such as GCM's. 

First let's review a popular approach introduced by \cite{wilks:05} that is strongly advocated in \cite{amp:13}. The key idea of this approach is to use a finite difference approximation to construct a time series that represents the model error (residual) when model \eqref{lor96_2} is used in place of the unknown dynamics in \eqref{lor96}. In this particular example, this residual time series can be obtained as follows:
\BEA
r_i(t) \approx \left( \frac{x_i(t+\delta t)-x_i(t)}{\delta t} \right) - x_{i-1}(t)(x_{i+1}(t)-x_{i-2}(t)) - x_i(t) + F,
\EEA
where we use a very short time step $\delta t=0.005$ relative to the Lorenz 3 days decaying time scale (or 0.2 model time unit) following \cite{amp:13}. Obviously, this is a poor approximation when the data is noisy or sparse in time; also, it requires knowing all components of $x_i(t)$. Then they apply a standard least squares method\index{least squares method} to fit this time series to a polynomial equation, such as,
\BEA
r_i =  -\zeta - \alpha x_i - \beta x_i^2 - \gamma x_i^3 + \tilde{r}_i,\label{cubicmodel}
\EEA
where the residuals $\tilde{r}_i$ from the polynomial fitting are subsequently fitted again to an AR(1) model, $\tilde r_i(t) = \phi \tilde r_i(t-\delta t) + \hat{\sigma}(1-\phi^2) \dot{W}_i(t)$, where $\dot{W}_i(t) \sim\mathcal{N}(0,t)$ are standard i.i.d white noises.
We should mention that this multiple regression fitting approach has been generalized and used to infer parameters of nonlinear multilevel regression models \index{multilevel regression models} \cite{kkg:05,kkg:06} where the model error estimators are chosen to linearly depend on $x$. Repeating this multiple regression procedure, we reproduce the parameters in \cite{amp:13}, which are $\zeta = -0.198$, $\alpha=0.575$, $\beta=-0.0055$, $\gamma=-0.000223$, $\phi = .993$, $\hat \sigma = 2.12$.\footnote{Here negative signs are used in \eqref{cubicmodel} for consistent notations throughout this note. In \cite{bh:14}, they presented the same results without negative signs.}   Let's denote this stochastic parametric model as the \emph{Cubic+AR(1)} reduced model. We found that this model is not useful at all for data assimilation when observations of $x_i$ are spatially sparse (a total $M=4$ observations resulting from observing at every other grid point of $x_i$), the filtered solutions with this model diverges catastrophically (the average RMSE goes to numerical infinity).  Now let's repeat the same fitting procedure on a simpler model error estimator, enforcing $\zeta=\beta=\gamma=0$ in \eqref{cubicmodel} and $\phi=0$ such that $\tilde{r}_i(t) = \hat{\sigma}\dot{W}_i(t)$; essentially, we want to fit a linear damping (where we hope that $\alpha>0$) and a white noise; in this case we obtain $\alpha=0.481$ and $\hat \sigma=2.19$ and let's call the resulting model as the \emph{offline} model. 

Alternatively, let's fit these two parameters adaptively or \emph{online}. Technically, we employ the state-augmentation approach to obtain $\alpha$ with the following model,
\begin{equation}\label{reducedmodel}
\begin{split}
\frac{d\tilde x_i}{dt}&= \tilde x_{i-1}(\tilde x_{i+1}-\tilde x_{i-2}) -\tilde x_i + F + \Big[- \alpha x_i(t) + \hat{\sigma} \dot{W}_i(t)\Big],\\
\frac{d\alpha}{dt} &= 0,
\end{split}
\end{equation}
and simultaneously implement the adaptive covariance estimation method discussed in Section~\ref{sec22} to obtain $\hat \sigma$ and the observation error covariance $R$. In \eqref{reducedmodel}, the terms in the square bracket in the dynamical equations for $x_i$ are the estimator for \eqref{edyn}. The implementation detail of this parameter estimation method is described in Appendix E of \cite{bh:14}. It is worth mentioning that the same strategy (fitting method) has been applied to parameterize the physics constrained models in \cite{mh:13,hmm:14} and to parameterize the Markovian models for the memory and noise terms in \eqref{edyn} as effective reduced models for Fourier modes of the Nonlinear Schr\"odinger equation \cite{hl:15}.   

For comparison, we also include the perfect model experiment with the full model in \eqref{lor96} which runs ETKF with an ensemble of size 528, doubling the total state variables $N(J+1)$, whereas the reduced model only uses ensemble of size 18, doubling the dimension of the augmented slow variables and one parameter $\alpha$, $(J+1)$. In Figure \ref{2layerL96comp} we compare the performance on the filtering experiment in the presence of model error for different observation time intervals. We see that the offline model gives worse performance relative to the observation in terms of RMSE.  On the other hand, the online model produces filtered solutions with RMSEs that are relatively close to those of the full model. We should point out that while the full filter and the offline method runs ETKF with known observation error covariance matrix $R$, the online method directly estimates $R$. Moreover, while the offline method requires a training data set of $x_i$ to estimate the parameters $\alpha$ and $\hat \sigma$, the online model uses only noisy sparse observations $v_i$ to estimate these same parameters on-the-fly. We show the estimated parameters $\alpha$ and $\hat \sigma$ to be compared with those from the offline estimates. Relative to the offline estimates, the online method produces smaller damping coefficient $\alpha$ (which means the online model retains more memory) and smaller noise amplitude $\hat \sigma$ (which implies that it is more accurate). The slight deterioration of the full model for long observation time relative to the online method could be due to the stiffness of the full model. 

\begin{figure}
\centering
\includegraphics[width=.45\textwidth]{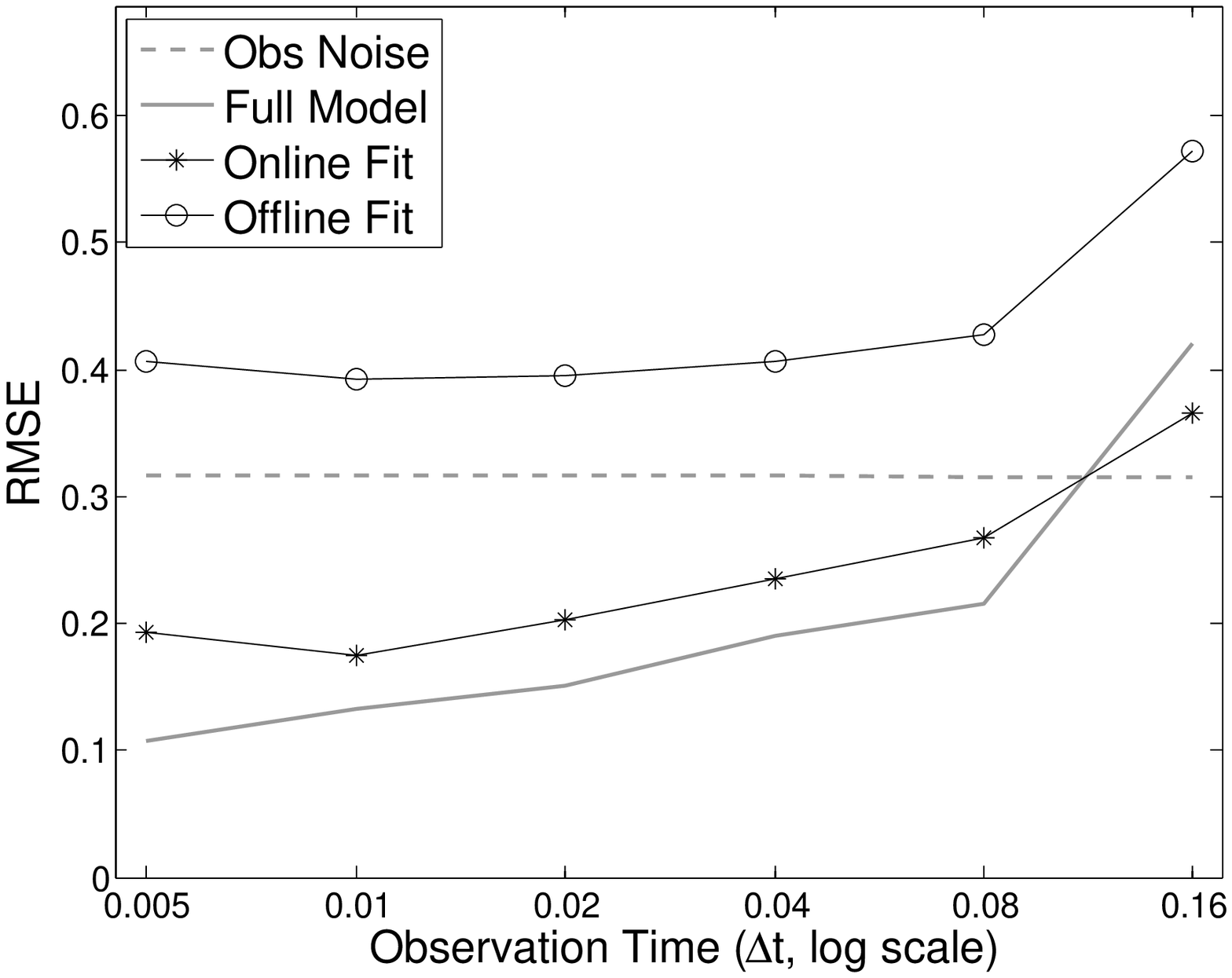}
\includegraphics[width=0.45\textwidth]{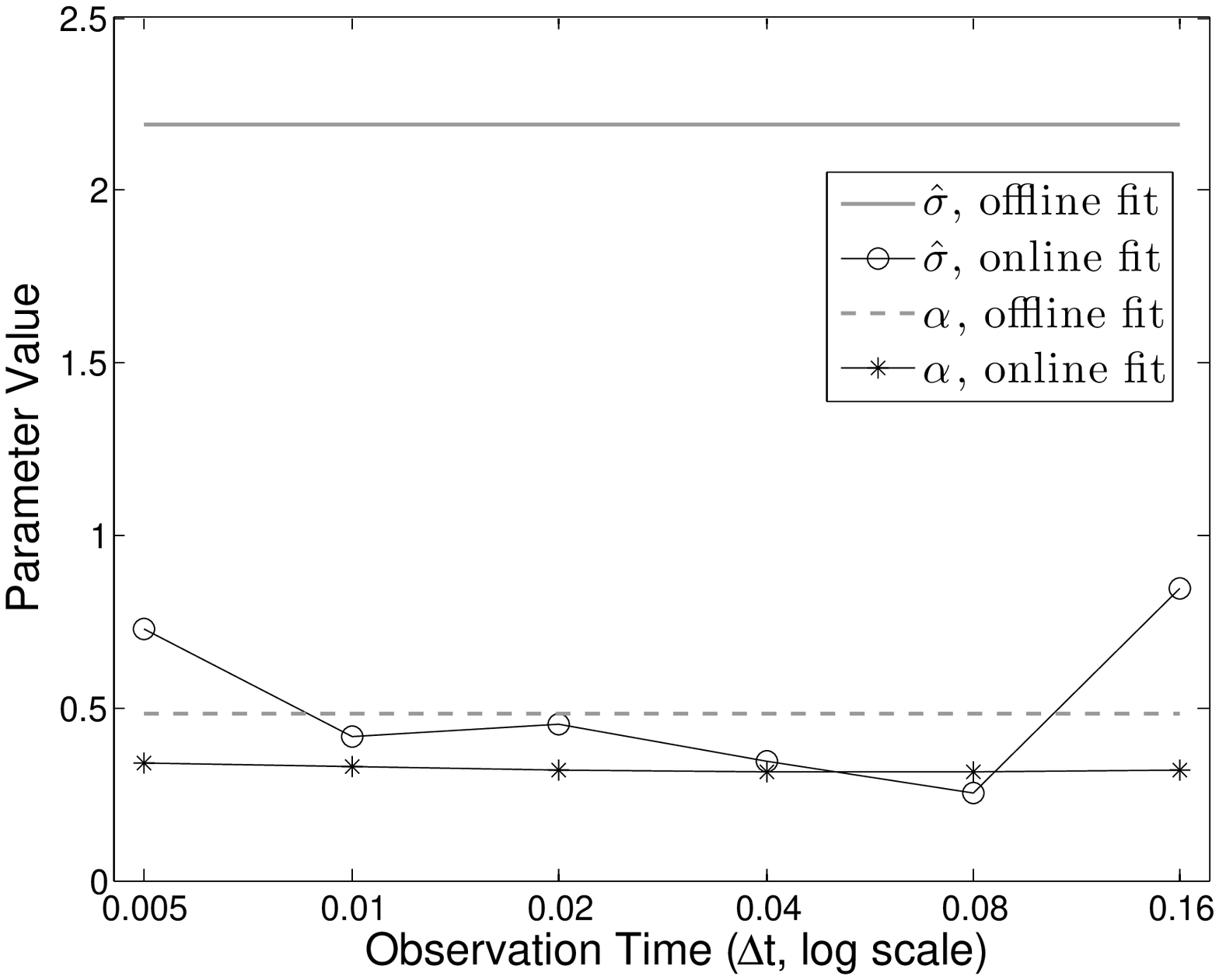}
\caption{\label{2layerL96comp} Filter performance measured in terms of root mean squared errors as functions of observation time interval (left). The full model filter uses \eqref{lor96}, the same model used to generate the data.  The \emph{Cubic+AR(1)} model is not shown since the filtered diverged. In the right panel, we compare the $\alpha$ and $\hat \sigma$ parameters from the online and offline estimation techniques.  The \emph{Offline Fit} curves use parameters $\alpha = 0.481$ and $\hat \sigma = 2.19$ estimated using the technique of \cite{amp:13}.}
\end{figure}

In Figure \ref{pdfcorr} we compare the equilibrium marginal density and the correlation function of $x_i$ from the online and offline models to those of the slow variables of the full model.  In this regime, both the equilibrium density and the correlation function from the online model agrees with those from the full model over a very long time (note that $4$ model time units corresponds to $800$ integration steps for the reduced model).  In contrast, the offline model and even the \emph{Cubic+AR(1)} model advocated in \cite{amp:13} showed some deviations, notably underestimating the variance and overestimating the lag correlations at the later times. Since the online model gives good filter performance and also closely matches the equilibrium statistics of the full model, we conclude that for this specific example,
the model error estimator can be modeled by a linear damping and a white noise in this regime. Furthermore, the online parameter estimation scheme is a natural way to infer the parameters in this stochastic model. The problem with the linear regression based estimation scheme of \cite{amp:13} is that the deterministic parameter, $\alpha$, and diffusion amplitude, $\hat{\sigma}$, in the stochastic parameterization model in \eqref{reducedmodel} are estimated separately. So, when a parameter in \eqref{reducedmodel} is independently perturbed, the nonlinear feedback of this perturbation is not appropriately accounted in the filtered estimates. In contrast, the online method constantly accounts for the nonlinear feedback of the perturbed parameters through the adaptive estimation strategy. 

\begin{figure}
\centering
\includegraphics[width=.45\textwidth]{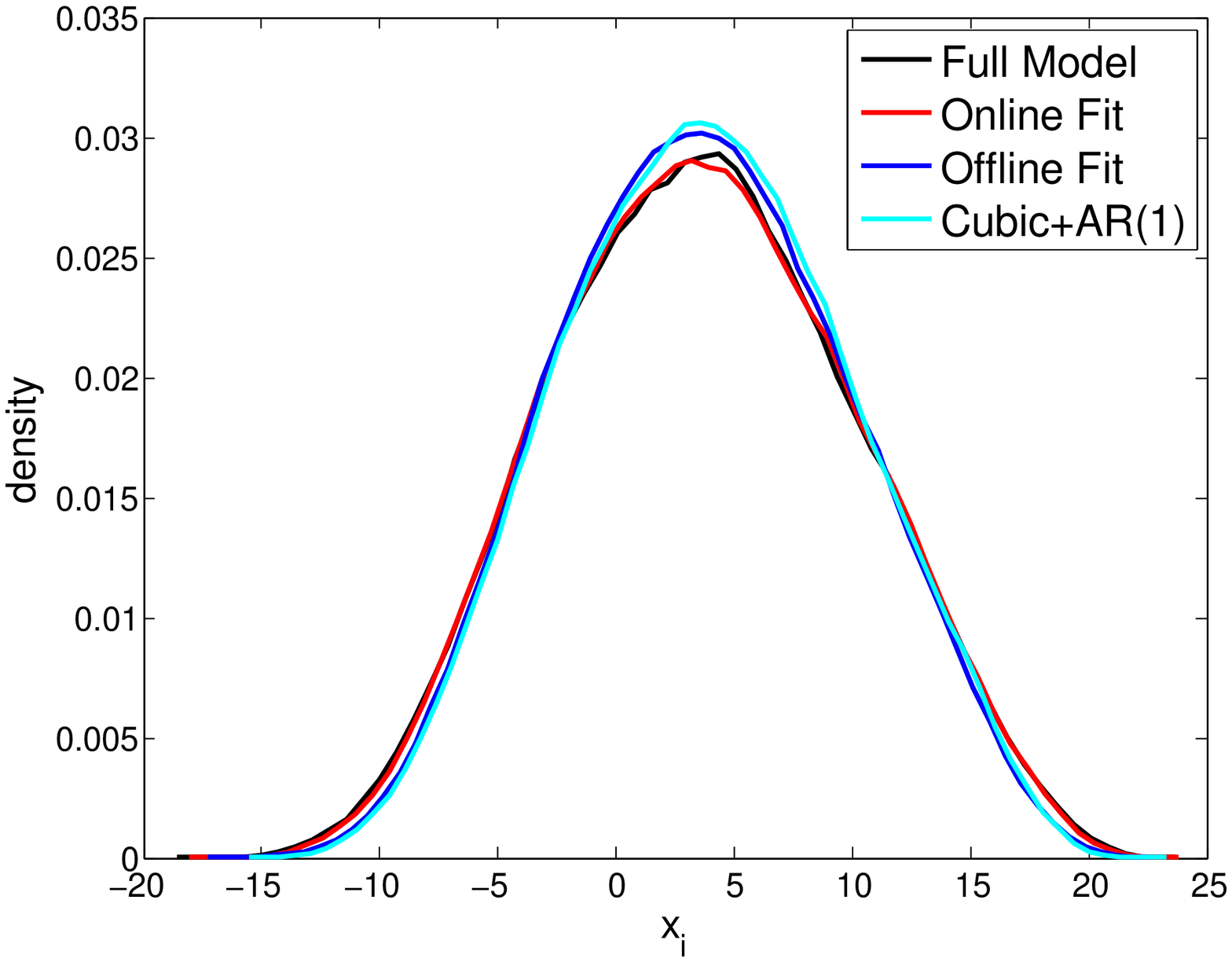}
\includegraphics[width=0.435\textwidth]{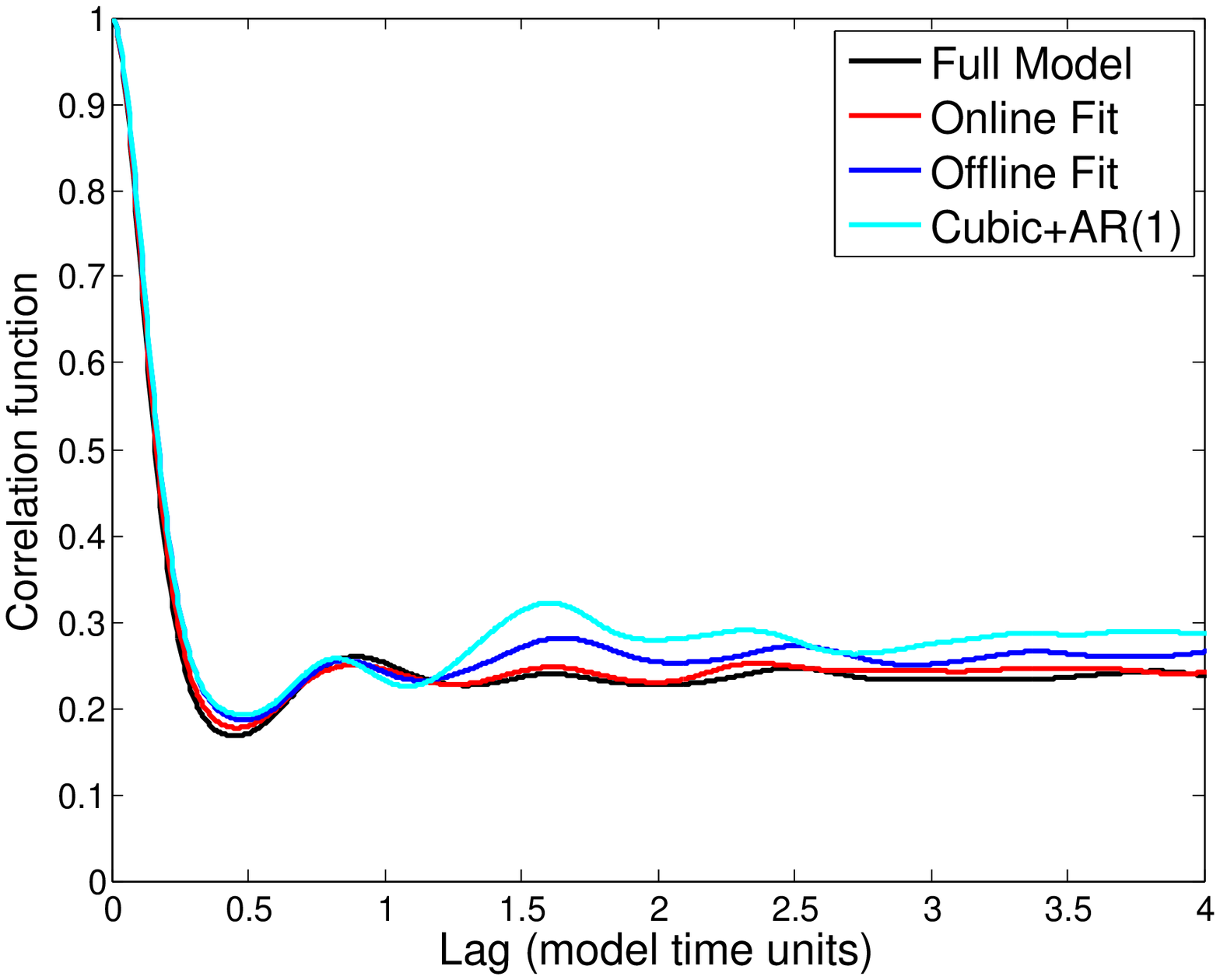}
\caption{\label{pdfcorr} Climatological forecast performance is compared in terms of the invariant measure of the slow variables shown as a probability density function (left) and the autocorrelation as a function of lag steps of length $0.005$ (right). Each curve is computed from a long free run with $1.6\times 10^7$ data points.}
\end{figure}
\end{example}
  
While the online parameterization methods discussed in the example above are basic recursive methods that estimate mean and covariance statistics (see Sections~\ref{sec21} and \ref{sec22}), it is worth mentioning that carefully designed offline minimization schemes have also been proposed to optimize some functionals such as the information theoretic \index{information theoretic}criteria \cite{mg:11a,mg:11}. We should note that such schemes are typically numerically expensive but the parameters are estimated offline (only once) and this approach has been shown to be useful in some applications \cite{cm:15}. 

A simple fact that one should notice from example 3 is that the online fitting through the filtering procedure produces statistics of the joint state-parameters $(x,\alpha,\hat\sigma)$ with respect to the posterior density, $p(x,\alpha,\hat \sigma |v)$. This suggests that one should consider investigating model error in data assimilation with a posterior distribution formulation rather than with the prior distribution formulation as in Section~\ref{sec1}. A natural question one will ask is: Why does the strategy in the example above work? In other words, why do the chosen stochastic model in \eqref{reducedmodel} simultaneously produce accurate filtering and climatological statistical predictions? Is this just a coincidence or can we justify this finding mathematically? We will address these questions in the next section.

\section{A linear theory for filtering with model error from unresolved scales}\label{sec4}

In this section, we review the theoretical result established in \cite{bh:14} which is based on analyzing model error from unresolved scales with a posterior distribution formulation.\index{posterior distribution formulation} We will compare it to that of a standard treatment from the prior distribution formulation \cite{gh:13}. Since our plan is to mathematically answer all of the questions above, we will consider a special setting relative to \eqref{slow}-\eqref{fast} to understand the issue. In particular, we consider stochastic dynamical systems of the following form,
\BEA
dx &=& f(x,y;\theta)\,dt + \sigma_x(x,y;\theta)dW_x, \label{slow2}\\
dy &=& \frac{1}{\epsilon}g(x,y;\theta)\,dt+\frac{\sigma_y(x,y;\theta)}{\sqrt{\epsilon}}dW_y,\label{fast2}
\EEA 
where $y$ is assumed to evolve on a faster time scale relative to $x$ and the scale gap is characterized by the parameter $\epsilon$. To facilitate the analysis in the simplest setting, we consider continuous-time observations,
\BEA
dz = x\,dt + \sqrt{R}dV,\quad R>0, \label{linobs}
\EEA
where $dW_x, dW_y, dV$ are i.i.d.~Wiener processes and $\theta$ denotes the true model parameters. We should note that while the analytical derivation was performed with continuous-time filter, the numerical verification in all of the examples below will be based on discrete-time filtering with large observation times. 

The main result from \cite{bh:14} loosely states that: \emph{There exists a reduced model that involves only $\tilde x$ of the following form:\index{reduced models}
\BEA
d\tilde{x} = \tilde{f}(\tilde{x},\Theta)\,dt + \tilde{\sigma}_{\tilde{x}}(\tilde{x},\Theta)dW,\label{averaged}
\EEA
where $\Theta$ depends on $\epsilon$ and $\theta$, such that the filter mean and covariance estimates resulting from the reduced filter in \eqref{averaged}, \eqref{linobs} are close to the corresponding posterior statistics obtained from the true filter in \eqref{slow2}-\eqref{fast2}, \eqref{linobs}. To clarify, the statistics of the reduced filter are defined with respect to conditional density $p(\tilde{x}|z)$ while the statistics of the full filter are defined with respect to conditional density $p(x,y|z)$. For the linear and Gaussian case, the resulting reduced model in \eqref{averaged} can be specified uniquely. The same unique set of parameters can also be found by matching equilibrium statistics of \eqref{slow2} and \eqref{averaged}. In other words, a consistent reduced model that simultaneously gives optimal filtering as well as accurate climatological prediction exists and is unique in Gaussian and linear setting.
}

While this result supports the finding in example~3, that is, such a consistent reduced model exists, this theory does not provide a general way to find the reduced model for every problem. On the other hand, the results in example~3 suggested that even if the correct ansatz is given (i.e., damping and white noise in this case), a natural method to obtain these parameters should be based on a filtering procedure that gives conditional estimates. For the linear and Gaussian setting, offline fit on the second order statistics are sufficient because other than being sufficient statistics, the covariance statistics are closed, meaning they do not depend on higher-order statistics as opposed to nonlinear problem as shown in example~1. We believe that this is the key factor that explains why the not so carefully designed offline fitting method, such as \cite{wilks:05,amp:13}, tends to produce inaccurate estimates even if the same parametric form (damping and white noise) is used.

Rather than re-deriving this result (as shown in \cite{bh:14}), we use two examples below to find the connection of this result with the discussions in the previous sections. The first example is the linear example studied in \cite{gh:13,bh:14} and the second one is a nonlinear problem introduced in \cite{ghm:10a,ghm:10b}. With these simple examples, we hope to elucidate the importance of posterior distribution formulation over the prior distribution formulation in accounting for model error in data assimilation of multiscale dynamical systems. Second, we want to emphasize that while the stochastic parameterization is a powerful tool that implicitly accounts for all nontrivial statistics of model error, there are still many remaining challenges in lifting this idea to solve general problems.  

\begin{example}[Example 4:]
Consider filtering a partially observed two-scale linear system of stochastic differential equations \cite{gh:13},
\BEA
dx &=& (a_{11} x + a_{12} y)\, dt + \sigma_x\, dW_x, \label{lslow}\\ 
dy &=& \frac{1}{\epsilon}(a_{21} x + a_{22} y)\, dt + \frac{\sigma_y}{\sqrt{\epsilon}}\, dW_y. \label{lfast}
\EEA
Here, $W_x, W_y$ are independent Wiener processes, the parameter $\epsilon$ characterizes the time scale gap between the variables $x\in \mathbb{R}$ and $y\in \mathbb{R}$. We assume throughout that $\sigma_x,\sigma_y \neq 0$ and that the eigenvalues of the matrix, 
\begin{align}
A = \begin{pmatrix}
       a_{11} & a_{12} \\ \frac{1}{\epsilon}a_{21} & \frac{1}{\epsilon}a_{22}
      \end{pmatrix}, 
\nonumber
\end{align}
are strictly negative, to assure the existence of a unique joint invariant density $\rho_\infty(x,y)$. Furthermore we require $\tilde{a}=a_{11}-a_{12}a_{22}^{-1}a_{21}<0$ to assure that the leading order slow dynamics,
\begin{align}
d\tilde{x} = \tilde a \tilde{x}\, dt + \sigma_x\, dW_x,\label{rsf}
\end{align}
supports an invariant density. It is well known that solutions of the one-dimensional SDE in \eqref{rsf} converge to solutions, $x^\epsilon(t)$, of \eqref{lslow} pathwise up to finite time, assuming $\epsilon\rightarrow 0$.  The convergence rate is on the order of $\epsilon$ (see e.g.,\cite{ps:00} for detail). Relating to \eqref{GLE} and the discussion preceding to \eqref{edyn}, one can think of $\bar{f}(\tilde{x}) = \tilde a \tilde{x} + \sigma_x\, \dot{W}_x$ as a result of the following projection $\bar{f}(\tilde x) = \lim_{\epsilon\rightarrow 0} \mathbb{E}[a_{11}x+a_{12}y+\sigma_x\dot{W}_x |\tilde{x}]$, where the expectation is taken with respect to the invariant density $p_\infty(y|x)$. Here, we use the physicist notation, for white noise $\dot{W}_x \equiv dW_x/dt$ to simplify the notation.

\noindent {\bf Reduced Stochastic Filter (RSF)}: \index{RSF (Reduced Stochastic Filter)}
Consider \eqref{rsf} as a prior model to assimilate noisy observations,
\begin{align}
v_m = x(t_m) + \varepsilon^o_m, \quad \varepsilon^o_m\sim\mathcal{N}(0,R),\label{obsmodel1}
\end{align}
of the slow variable $x$ at discrete time step $t_m$ with constant observation time interval $\Delta t=t_{m+1}-t_m$. Connecting to the discussion before \eqref{edyn}, this approach essentially offers no treatment on model error, that is, $e=0$ since $\bar{f}(\tilde{x}) = \tilde a \tilde{x}+\sigma_x\, \dot{W}_x$.
Since this example is linear, the optimal solutions can be obtained by the Kalman filter formula, in the sense that the solutions minimize the posterior error variance \cite{kalman:61}. In discrete form, the prior mean and error covariance estimates \cite{gardiner:97,mh:12} are given by
\begin{align}
\bar{\tilde{x}}^b_{m} &= F \bar{\tilde{x}}^a_{m-1},\nonumber\\
\tilde{P}^b_{m} &= F\tilde{P}^a_{m-1}F^\top + Q,\nonumber
\end{align}
where $F=e^{\tilde{a}\Delta t}$ and $Q = \frac{\sigma_x^2}{-2\tilde{a}}(1-e^{2\tilde{a}\Delta t})$. We should emphasize that this $Q$ is not associated with statistics of model error. This $Q$ is the variance of the stochastic forcing in the reduced model in \eqref{rsf}. The posterior mean and covariance update are given by,
\begin{align}
\bar{\tilde{x}}^a_{m} & = \bar{\tilde{x}}^b_{m} + K_m(v_m - \bar{\tilde{x}}^b_m),\nonumber\\ 
\tilde{P}^a_m &= (1-K_m)\tilde{P}^b_m, \nonumber\\
K_m &= \tilde{P}^b_m(\tilde{P}^b_m+R)^{-1}.\nonumber
\end{align}
We will refer to this filtering scheme as the reduced stochastic filter (RSF) as in \cite{gh:13}. It has been shown that 
the posterior filtered estimates of such a reduced stochastic filter converge to the true filtered solutions, with a convergence rate of $\sqrt{\epsilon}$ for general nonlinear filtering problems, see \cite{imkeller:13}.   

Now we discuss results from a numerical simulation with $a_{11}=a_{21}=a_{22}=-1, a_{12}=1$, $\sigma_x^2=\sigma_y^2=2$, $\Delta t=1$, and $R=50\%Var(x)$ and compare them with the {\bf true filtered solutions}, obtained with the perfect prior model in \eqref{lslow}-\eqref{lfast}. In Figure~\ref{fig1}, we show the filter accuracy (left panel), quantified by the Mean-Square-Error (MSE) \index{MSE (Mean-Square-Error)} between the posterior mean state estimate, $\bar{\tilde{x}}^a_m$,  and the truth, $x_m$, and the asymptotic error covariance estimate (right panel) of the posterior mean estimate, $\bar{\tilde{x}}^a_m$, as functions of scale gap $\epsilon$.
Note that the asymptotic posterior error covariance estimate is constant for this linear problem after $m=10,000$ iterations. Notice also that when $\epsilon\ll 1$ is small ($x$ is much slower than $y$), the MSEs are almost identical to those of the true filter. For moderate scale gap with larger $\epsilon$, notice that the filter accuracy degrades (with higher MSE) and the true prior error covariance $P^b_m$ is significantly underestimated (see the solid line with circles in Figure~\ref{fig1}). 

\begin{figure*}
\centering
\includegraphics[width=0.45\textwidth]{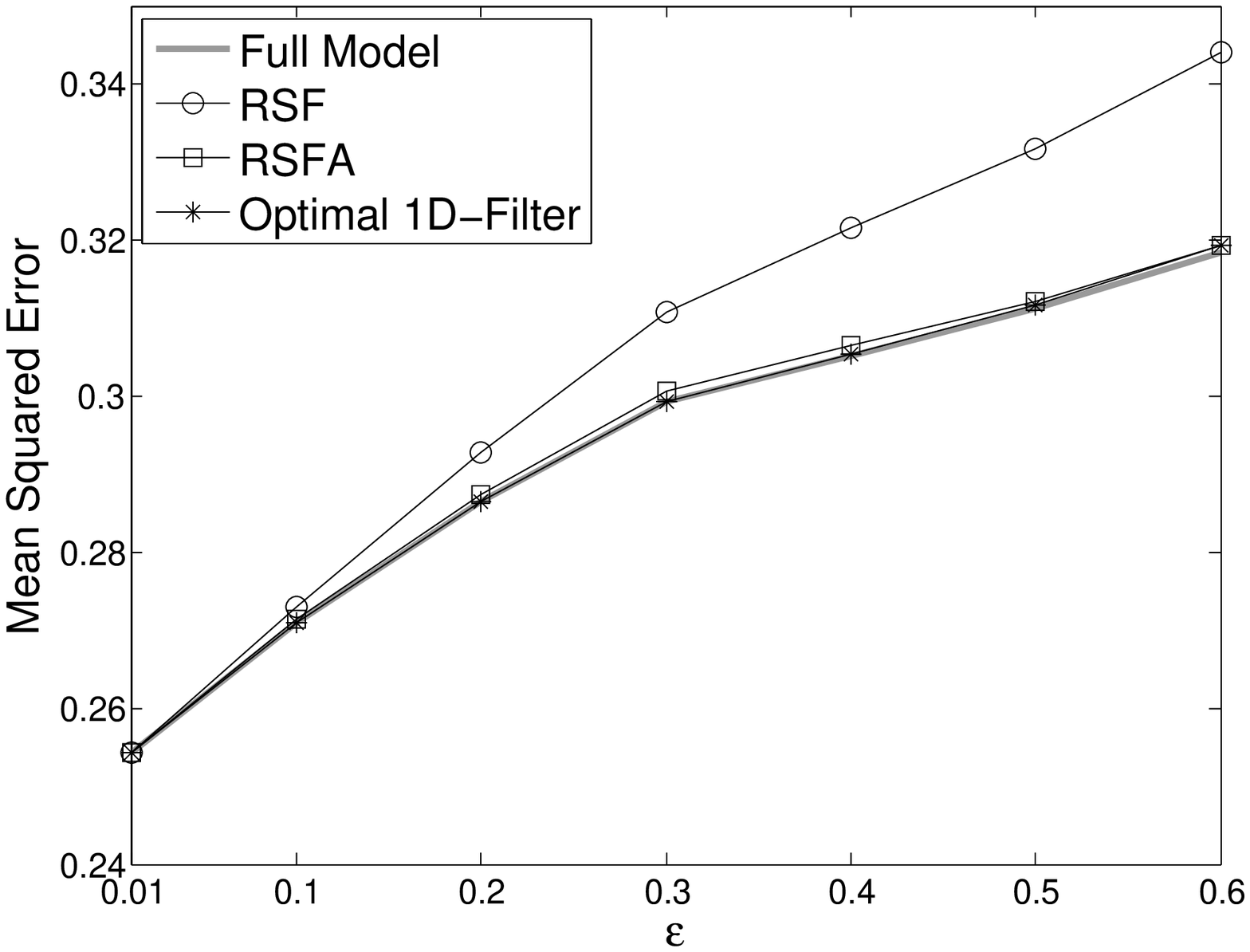}
 \includegraphics[width=0.45\textwidth]{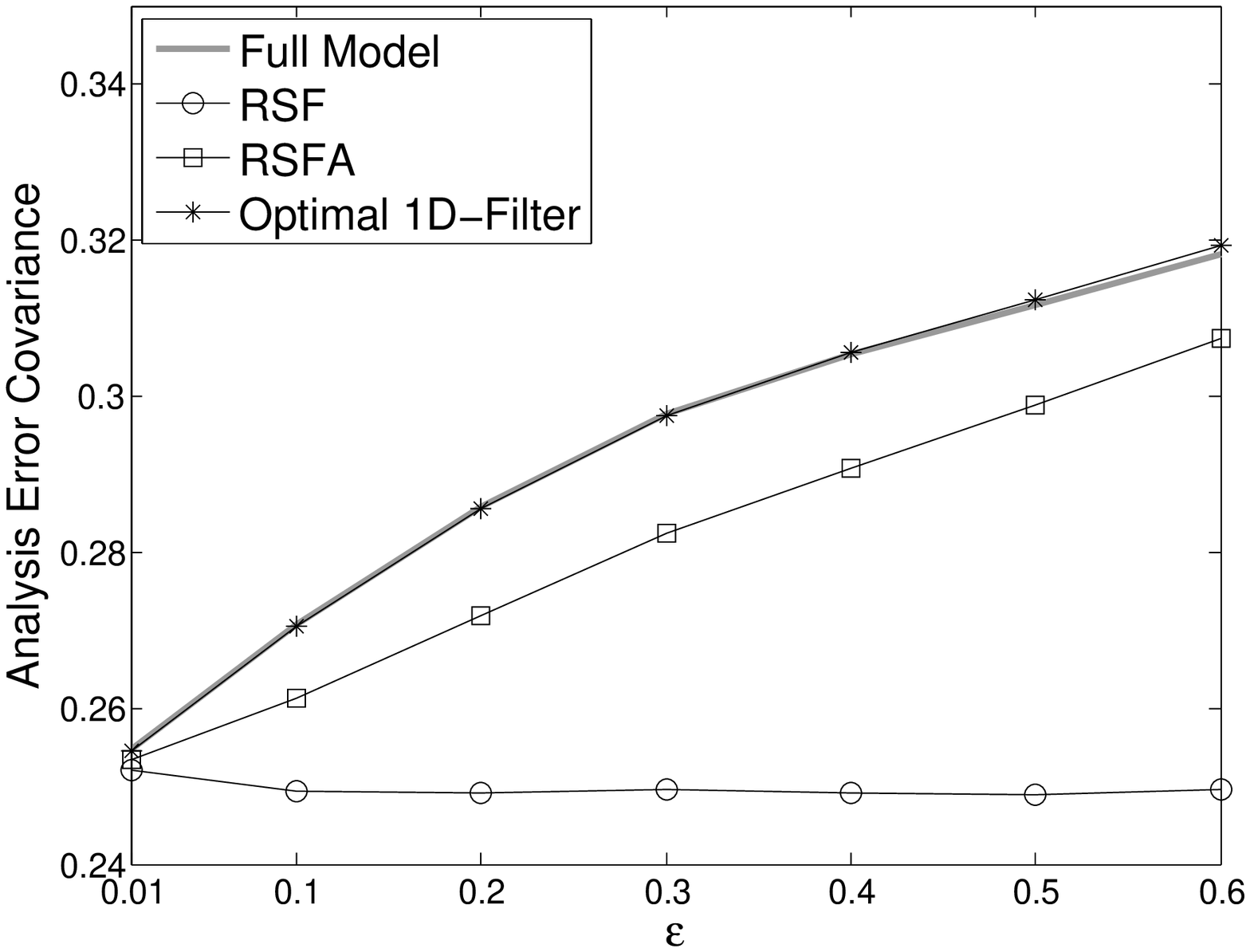}
\caption{Average mean square error (left panel) and the asymptotic posterior error covariance estimate (right panel) as functions of scale gap $\epsilon$ for filtering the linear problem in \eqref{lslow}-\eqref{lfast}.}
\label{fig1} 
\end{figure*}

\noindent{\bf RSF with an additive noise correction (RSFA):}\index{RSFA (Reduced Stochastic Filter with Additive noise correction)} Here, we will use a prior distribution formulation to treat model error. In particular, we just apply an asymptotic expansion on the prior model, ignoring the availability of noisy observations. Let's rewrite the fast equation in \eqref{lfast} as follows,
\BEA
y\, dt = \frac{a_{21}}{a_{22}} x\,dt - \sqrt{\epsilon}\sigma_y\frac{a_{12}}{a_{22}}\, dW_y + \mathcal{O}(\epsilon).
\EEA
Substitute this expression into the slow equation in \eqref{lslow}, we obtain:
\begin{align}
d\hat{x} = \tilde a \hat{x}\, dt + \sigma_x\, dW_x - \sqrt{\epsilon}\sigma_y\frac{a_{12}}{a_{22}}\,dW_y.\label{rsfa}
\end{align}
One can check \cite{gh:13} for a more concise formal asymptotic expansion; therein, they also showed that solutions of \eqref{rsfa} converge pathwise to solutions, $x^\epsilon(t)$, of \eqref{lslow} up to finite time, with convergence rate of order $\epsilon^2$. We will refer to the filtering strategy with the prior model in \eqref{rsfa} as the reduced stochastic filter with an additive noise correction (RSFA), following the notation in \cite{gh:13}. The additional additive noise correction in \eqref{rsfa} essentially inflates the prior covariance estimates in each filtering step, so RSFA is an analog of an additive covariance inflation method \cite{ott:04,whwst:08}.

To be consistent with the notations in \eqref{GLE}, \eqref{edyn}, and assuming that $\bar{f}(\tilde x)=\tilde a\tilde x+\sigma_x\dot{W}_x$, this reduced filter model is equivalent to modeling \eqref{edyn} with the following estimator, 
\begin{align}
d\hat{e} = \tilde{a}\hat{e}\,dt - \sqrt{\epsilon}\sigma_y\frac{a_{12}}{a_{22}}\,dW_y.\label{bhat1}
\end{align}
where $\hat{e}\equiv\hat{x}- \tilde{x}$; here, $\tilde{x}$ solves \eqref{rsf} and $\hat{x}$ solves \eqref{rsfa}. When $\hat{e}(0)=0$, the model error estimator is essentially an unbiased white noise process. Our numerical simulations suggest that while the filter accuracy is improved (notice in Figure~\ref{fig1} that the MSE are almost identical to those of the true filter), the true posterior error covariances, $P^b_m$, are still underestimated. 

\noindent {\bf Optimal Reduced Stochastic Filter:}\index{optimal reduced stochastic filter} Finally, let's discuss the model error estimator resulting from a posterior distribution formulation on the continuous-time filtering problem in \eqref{lslow}, \eqref{lfast}, \eqref{linobs}. In \cite{bh:14}, they rigorously proved that there exists a unique choice of estimator of model error, $e\equiv x-\tilde{x}$, such that the filtered solutions are optimal in the sense that both the mean and covariance estimates are as accurate as those of the true filter. The model error estimator for \eqref{edyn} satisfies the following dynamics,
\begin{align}
d\hat{e} = \tilde{a}\hat{e}\,dt  - \sqrt{\epsilon}\sigma_y\frac{a_{12}}{a_{22}}\,dW_y -\epsilon \hat{a}\tilde{a}(\hat{e} + \tilde{x})\,dt-\epsilon\sigma_x\hat{a}\,dW_x.\label{de}
\end{align}
where $\hat{a}\equiv a_{12}a_{21}/a_{22}^2$. Notice that the model error estimator is not just a Gaussian white noise, in this case it also depends on $\tilde x$. With this model error estimator, the reduced filter prior model is given by
\begin{align}
d\hat{x} &= d\tilde{x} + d\hat{e}\nonumber\\
&=\Big( \tilde a \tilde{x}\, dt + \sigma_x\, dW_x\Big) + \Big(\tilde{a}\hat{e}\,dt  - \sqrt{\epsilon}\sigma_y\frac{a_{12}}{a_{22}}\,dW_y -\epsilon \hat{a}\tilde{a}(\hat{e} + \tilde{x})\,dt-\epsilon\sigma_x\hat{a}\,dW_x\Big)\nonumber\\
&= \tilde{a}(1-\epsilon\hat{a})\hat{x}\,dt + \sigma_x (1-\epsilon\hat{a})\, dW_x - \sqrt{\epsilon}\sigma_y\frac{a_{12}}{a_{22}}\,dW_y,\label{optimalreducedmodel}
\end{align}
where $\hat{x}\equiv\tilde{x}+\hat{e}$. 

We numerically confirm the accuracy of both the mean and covariance estimates with this optimal reduced model in Figure~\ref{fig1}. We should also point out that this result was found by enforcing the linear optimality condition, $\mathbb{E}(e\cdot \hat{x})=0$ (which is satisfied when a filtered mean estimate is optimal \cite{oksendal:03}). With this choice of parameters, the reduced filtered solutions become consistent in the sense that the actual error covariance of the filtered mean estimate matches the filtered error covariance estimate, $\mathbb{E}[e^2]=\mathbb{E}[(x-\tilde{x})^2]+\mathcal{O}(\epsilon^2)$. Numerically, notice that the MSE (a numerical estimate for the actual error covariance estimate) and the posterior error covariance estimate in Figure~\ref{fig1} are very similar for only the true filter and the optimal one-dimensional filter. In this example, these are the only consistent filters.

As we pointed out before, the same reduced model in \eqref{optimalreducedmodel} can be determined by fitting the reduced filter model to the equilibrium covariance statistics and the correlation time of the underlying true signal that solves \eqref{lslow}-\eqref{lfast} for the slow variable $x$. As a consequence, the optimal reduced model in \eqref{optimalreducedmodel} produces, both, an optimal filtering and an optimal equilibrium statistical prediction. This is the linear theory established in \cite{bh:14}. In this linear and Gaussian setting, parameters of the reduced model can be obtained offline by fitting climatological statistics. 
\end{example}

While numerical example 3 suggests a possibility for this theory to hold for general nonlinear systems, the problem becomes much more difficult to analyze in a general setting. For general continuous-time nonlinear filtering problems, the true filtered solutions are characterized by conditional densities, which solve a stochastically forced partial differential equation known as the Kushner equation \cite{kushner:64}\index{Kushner equation} and solving the Kushner equations is nontrivial for general high-dimensional nonlinear problems. Rather than attempting to analyze this issue in a general setting since it may not necessarily give practical algorithms to tackle high-dimensional problems, we will use the next example to verify the linear theory above on a simple nonlinear test model.

\begin{example}[Example 5:]
Consider the nonlinear filtering problem \cite{gh:13} of noisy observations,
\begin{align}
v_m = x(t_m) + \varepsilon^o_m, \quad \varepsilon^o_m\sim\mathcal{N}(0,R),\label{obsmodel2}
\end{align}
where
\begin{align}
\frac{dx}{dt} &= -(\gtilde+\lh)x + \btilde + f(t) +\sx \dot{W}_x,  \nonumber \\
\frac{d\btilde}{dt} &= -\frac{\lb}{\epsilon}\btilde + \frac{\Sb}{\sqrt{\epsilon}} \dot{W}_b, \label{spekf} \\
\frac{d\gtilde}{dt} &= -\frac{\dg}{\epsilon} \gtilde  + \frac{\sg}{\sqrt{\epsilon}} \dot{W}_\gamma, \nonumber
\end{align}\index{SPEKF (Stochastic Parameterized Extended Kalman Filter)}
with $\lh = \gh-\mathi \omega$ and $\lb = \gb-\mathi \wb$. The model in \eqref{spekf} was introduced as a test model of a stochastic parameterization for filtering a turbulent mode in the presence of model error in \cite{ghm:10a,ghm:10b}. The solutions for the nonlinear filtering problem in \eqref{spekf}, \eqref{obsmodel2}, was called SPEKF, which stands for Stochastic Parameterized Extended Kalman Filter \cite{ghm:10a,ghm:10b,mhg:10,mh:12}. In particular, SPEKF posterior statistical solutions are obtained by applying Kalman update to the exactly solvable prior statistical solutions of \eqref{spekf}. We should point out that the SPEKF solutions are \emph{not} the true filtered solutions. For general continuous-time nonlinear filtering problems, the true filtered solutions are characterized by the conditional distribution $p(x_t,\btilde_t,\gtilde_t| z_\tau, 0\leq \tau\leq t)$, which solves the Kushner equation \cite{kushner:64}. It turns out that the posterior solutions of SPEKF for discrete observation time are the analog of the Gaussian closure on the first two-moments of this conditional distribution for the corresponding continuous-time filter \cite{bh:14}. In this sense, one can refer to SPEKF solutions as the best approximate solutions that are numerically attainable since the true filtered solutions are not accessible.   

The nonlinear system in \eqref{spekf} has many attractive features as a test model. First, it has exactly solvable statistical solutions which are non-Gaussian. Thus, it allows one to study non-Gaussian prior statistics conditional to the Gaussian posterior statistical solutions of the Kalman update and to verify uncertainty quantification methods \cite{bm:13}. 
Second, a recent study by \cite{bgm:12} suggests that the system in \eqref{spekf} can reproduce signals in various turbulent regimes such as intermittent instabilities in a turbulent energy transfer range and in a dissipative range as well as laminar dynamics. 

As in the linear example~4 above, the $\mathcal{O}(1)$ dynamics are given by the averaged dynamics, where the average is taken over the unique invariant density generated by the fast dynamics of $\btilde$ and $\gtilde$ \cite{gh:13}, which results in a linear SDE,\index{RSF (Reduced Stochastic Filter)}
\begin{align}
\frac{d\tilde{x}}{dt}  =  -\lh  \tilde{x} + f(t) + \sx \dot{W}_x.\label{e.nlav}
\end{align}
In the numerical simulation below, we will refer to the filtering scheme with the prior model in \eqref{e.nlav} as the Reduced Stochastic Filter (RSF). This approach essentially offers no model error treatment, assuming that the right-hand-terms in \eqref{e.nlav} is $\bar{f}(\tilde{x})$ (to be consistent with our previous notations in \eqref{edyn}).   In \cite{gh:13}, they defined a reduced stochastic filter with an additive noise correction (RSFA) given by the following model error estimator, \index{RSFA (Reduced Stochastic Filter with Additive noise correction)}
\begin{align}
\frac{d\hat{e}}{dt} = -\hat{\lambda}\hat{e} + \sqrt{\epsilon} \frac{\Sb}{\lb} \dot{W}_b.\label{bias3}
\end{align}
When initial model error is absent, $\hat{e}(0)=0$, this formulation essentially approximates the memory and noise terms in \eqref{edyn} with an unbiased white noise process. 

The posterior distribution formulation in \cite{bh:14}, suggested that the best one-dimensional reduced filtering (best in the sense that the errors in mean and covariance are of order-$\epsilon$ from the solutions of SPEKF) can be achieved with a damping and combined, additive and multiplicative, noise corrections,
\begin{align}
\frac{d\hat{e}}{dt} = -\hat{\lambda}\hat{e} + \sqrt{\epsilon}\Bigg( \frac{\Sb}{\sqrt{|\lb(\lb +\epsilon\lh )|^2}} \dot{W}_b - \frac{\sg}{\sqrt{\dg(\dg+\epsilon\gh )}}(\tilde{x}+\hat{e}) \circ \dot{W}_\gamma\Bigg),\label{bias4}
\end{align}
where the multiplicative noise term in \eqref{bias4} is Stratonovich. Notice that we refrain from calling the model estimator in \eqref{bias4} the optimal estimator since the optimal filtered solutions are not accessible unless  one can solve the Kushner equation for the full conditional distribution as we explained above. We will refer to the filtered solutions corresponding to model error estimator in \eqref{bias4} as the reduced SPEKF solutions. \index{reduced SPEKF}
 
Notice that when $\epsilon\gh \ll \dg$, the noise correction model in \eqref{bias4} can be approximated by,
\begin{align}
\frac{d\hat{e}}{dt} = -\hat{\lambda}\hat{e} + \sqrt{\epsilon} \frac{\Sb}{\lb} \dot{W}_b - \sqrt{\epsilon}\frac{\sg}{\dg}(\tilde{x}+\hat{e}) \circ \dot{W}_\gamma,\label{bias5}
\end{align}
which yields the reduced stochastic prior model RSFC, introduced in \cite{gh:13}. We should point out that the multiplicative noise in \cite{gh:13} is also in the Stratonovich sense. \index{RSFC (RSF with Combined additive \& multiplicative noises)}
Comparing the model error estimator in \eqref{bias4} and \eqref{bias5}, we notice that while it is possible to obtain the same parametric form (damping, additive and multiplicative noise forcings) by formulating through either the posterior and prior distribution, the resulting parameters in the two estimators are very different if condition $\epsilon\gh \ll \dg$ is not satisfied. In \cite{bh:14}, it was shown that for a set of parameters corresponding to dissipative range, in which the condition $\epsilon\gh \ll \dg$ is not satisfied, the resulting reduced model in \eqref{bias5} produces covariance statistics that are unstable. The main point we want to make is that the posterior distribution formulation provides more robust model error estimators.

Here, we only show the numerical results for the parameter set corresponding to the turbulent transfer energy range regime
\cite{bgm:12,gh:13}, $\epsilon = 1$, $\gh=1.2$, $\gb=0.5$, $\dg=20$, $\sx=0.5, \Sb=0.5, \sg=20$. In this regime, $x(t)$ exhibits frequent rapid transient instabilities, and $\gtilde$ decays faster than $u$, that is, $\epsilon\gh<\dg$, such that the RFSC in \eqref{bias5} is a good approximation of the reduced SPEKF with model error estimator \eqref{bias4}. The noisy observations in \eqref{obsmodel2} are sampled at every time interval $\Delta t=0.5$ (shorter than the decay time 0.833) and the noise variance is $R = 0.5Var(u)$. We will show the numerical results of three reduced filters, where the analyses are updated by the Kalman filter formula with prior models: (i) RSF in \eqref{e.nlav}, (ii) RSFA, accounting for model error with the stochastic model in \eqref{bias3}, and (iii) reduced SPEKF, accounting for model error with the stochastic model in \eqref{bias4}. We compare the estimates from these three filters with those from solutions of SPEKF in Figure~\ref{fig2}.

\begin{figure}
\centering
\includegraphics[width=0.7\textwidth]{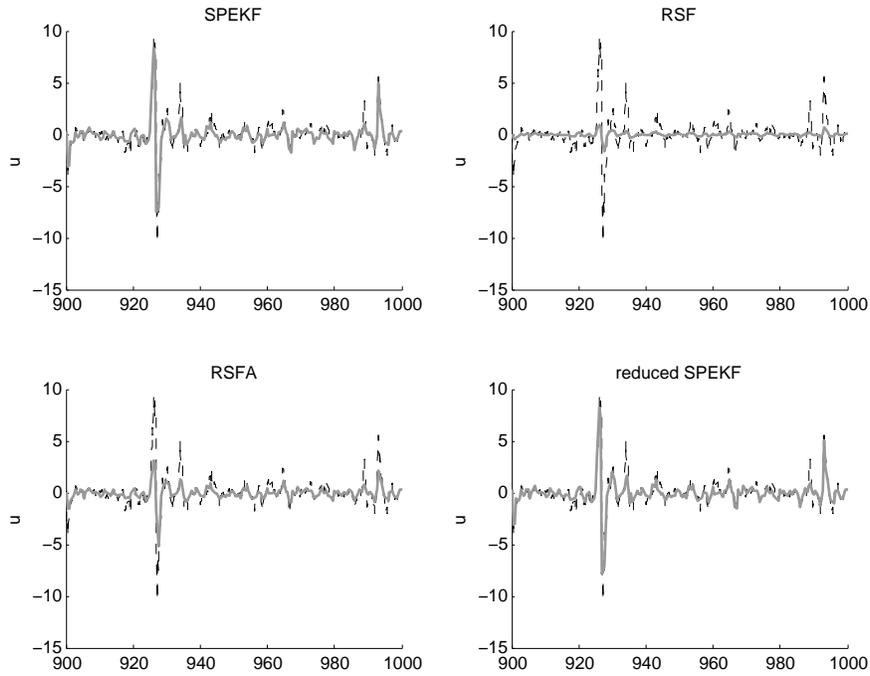}
\caption{Trajectory of the posterior mean estimates (in grey) compared to the truth (dashes). The average RMS errors of are 0.7730 (SPEKF), 1.5141 (RSF), 1.1356 (RSFA), 0.7861 (reduced SPEKF), and the observation error is $\sqrt{R}=1.19$ as a reference.}
\label{fig2} 
\end{figure}

\begin{figure*}
\centering
\includegraphics[width=0.6\textwidth]{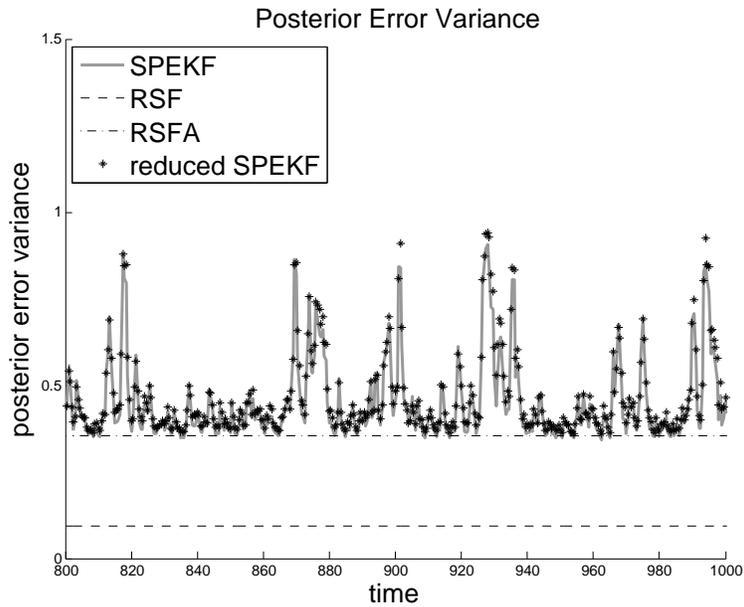}
\caption{Trajectory of the posterior covariance estimates corresponding to the filtered mean estimates in Figure~\ref{fig2}.}
\label{fig3} 
\end{figure*}

Notice that the reduced SPEKF, which accounts for model error with the combined additive and multiplicative noise in \eqref{bias4}, is the only one method which produces filtered solutions with accuracy that is comparable to that of SPEKF solutions (see Figure~\ref{fig2}); the average RMS errors (over 2000 iterations) are 0.7730 for the true filter, 0.7861 for the optimal filter, 1.1356 for RSFA, and 1.5141 for RSF. In Figure~\ref{fig3}, we show the corresponding posterior error covariance estimates from various reduced filters, $\tilde{P}^a_m$, compared to that from SPEKF, $P^a_m$ (in grey). Notice that RSF and RSFA significantly underestimate the posterior error covariances. The reduced SPEKF, on the other hand, tracks the covariance estimates from SPEKF, quite accurately. More comprehensive results based on various parameter regimes are shown in \cite{bh:14}; in there, they also showed that the same reduced model corrected with \eqref{bias4} produces accurate long term covariance solutions, with accuracy of order-$\epsilon$. 
\end{example}

These examples show the existence of a consistent reduced model based on a posterior distribution formulation and we verified that the theory is extendable on a special analytically tractable low-dimensional nonlinear model. While this theoretical result provides a firm understanding and more reason to use stochastic parameterization beyond white noise modeling, in general, the problems are still difficult to analyze. We suspect that the resulting model error estimator is much more complicated than \eqref{de} or \eqref{bias4}, with nontrivial dependence on $\tilde{x}$. For practical implementation, besides designing effective parameterization schemes as suggested in example~3, a more important wide open problem is to find a recipe to decide which parametric form is the adequate model in the sense that it is stable and it produces consistent equilibrium distribution and optimal filtering. Since choosing appropriate parametric models is difficult, we will discuss an alternative approach that does not use parametric modeling approach in the next section. 

\section{A nonparametric approach}\label{sec5}\index{nonparametric modeling}

From the discussion above, we showed that while parametric models can be used to learn high-dimensional systems from small data sets, their rigid parametric structure makes them vulnerable to model error. In this section, we overview a recently developed nonparametric modeling approach \cite{bgh:15,bh:15physd} and discuss a semiparametric framework, using the nonparametric approach as a model error estimator for dynamically evolving parameters of physics based parametric models. This framework was designed to avoid the two practical issues in stochastic parameterization, choosing appropriate parametric models and estimating the corresponding parameters. 

\subsection{The diffusion forecasting model}\label{sec51}\index{diffusion forecasting model}

Consider diffusion processes $\theta(t)$ that satisfy,
\begin{align}\label{SDE} d\theta = a(\theta)\, dt + b(\theta)\, dW_t, \end{align}
for a generic initial condition such that \eqref{SDE} is ergodic on a Riemannian manifold\index{Riemannian manifold} $\mathcal{M}\subset\mathbb{R}^n$. We assume that the distribution of $\theta$ can be characterized by a time-dependent density function, $p(\theta,t)$, that solves the Fokker-Planck equation,\index{Fokker-Planck equation}
\begin{align}
\frac{\partial p}{\partial t} = \mathcal{L}^*p = \nabla\cdot (-ap + \frac{1}{2}\nabla (bb^\top)p)  , \quad p(\theta,t) = p_t(\theta),\label{FokkerPlanck}
\end{align}
where $\mathcal{L}^*$ denotes the linear Fokker-Planck operator; here the differential operators are defined on $\mathcal{M}$ with respect to the Riemannian metric\index{Riemannian metric} inherited from the ambient space $\mathbb{R}^n$. We note that the equilibrium distribution, $\peq(\theta)$, of the underlying dynamics \eqref{SDE} satisfies, $\mathcal{L}^*\peq = 0$. In \eqref{SDE} and \eqref{FokkerPlanck}, $a(\theta)$ is a vector field which represents the deterministic part of the dynamics of $\theta$, and $b(\theta)$ is a diffusion tensor which determines the covariance structure of the stochastic forcing, $W_t$, which is a standard Brownian process on the manifold $\mathcal{M}$. 

Given a time series $\theta_i = \theta(t_i)$ sampled at discrete times $\{t_i\}_{i=1}^{N}$ we are interested in constructing a forecasting model so that given an initial density $p(\theta,t)$ at time $t$ we can estimate the density $p(\theta,t+\tau)$ at time $t+\tau$, where $\tau>0$. The key idea of the \emph{diffusion forecast} introduced in \cite{bgh:15} is to project the forecasting problem \eqref{FokkerPlanck} onto a basis of smooth real-valued functions $\{\varphi_j(\theta)\}$ defined on the manifold $\mathcal{M}$. Particularly, they chose $\{\varphi_j(\theta)\}$ to be the eigenfunctions of an elliptic operator $\hat{\cal L}$, corresponding to the generator \index{generator}of the gradient flows with isotropic diffusion\index{gradient flows with isotropic diffusion},
\BEA
d\theta = - \nabla U(\theta)\,dt + \sqrt{2} dW_t,\label{gradflow}
\EEA	
of the following potential function, $U(\theta) = -\log (\peq(\theta))$. The main motivation to choose these basis functions is that they are obtainable via the diffusion maps algorithm\index{diffusion maps algorithm} for data lying on compact manifold \cite{cl:06} and non-compact manifold \cite{bh:15vb}. In this presentation, we will use the algorithm in \cite{bh:15vb} since we are interested in the case where the sampling measure of $\theta$ are arbitrarily small and positive, which means $\mathcal{M}$ is non-compact. A short summary of the theory in \cite{bh:15vb} is given in the Appendix A of \cite{bgh:15} and the detail pseudo algorithm for constructing ${\varphi(\theta)}$ is presented in \cite{bh:15physd}. A second less obvious reason (which we will clarify below) is that this choice of eigenfunctions also minimizes the Dirichlet energy norm (with respect to $\peq(\theta)$) which turns out to minimize the stochastic error term in approximating the semigroup solutions of the adjoint of the Fokker-Planck operator of the general diffusion processes in \eqref{SDE}. Given all these facts, it is not difficult to show that the solutions of \eqref{FokkerPlanck} can be formally rewritten as follows (see \cite{bgh:15,bh:15physd} for details):   
\BEA
p(\theta,t+\tau) &=& \sum_j c_j (t+\tau) \varphi_j(\theta)\peq(\theta), \label{forecast}
\EEA
where coefficients $c_j(t+\tau) =\sum_l \langle \varphi_l,e^{\tau\mathcal{L}}\varphi_j\rangle_{\peq} c_l(t)$ and $c_l(t)=\langle p_t, \varphi_l\rangle$ will be numerically realized by Monte-Carlo approximations. Numerically, this approach can be interpreted as solving the linear Fokker-Planck equation with a spectral method\index{spectral method} in which the basis functions are estimated from data set $\theta_i$ without knowing $\mathcal{M}$. Practically, the diffusion maps algorithm will produce eigenvectors $\vec\varphi_j$ whose $i$-th component is an estimate of the eigenfunction $\varphi_j(\theta_i)$ evaluated at data set $\theta_i$. This is in contrast to the standard spectral methods \cite{trefethen:00} which impose a certain set of basis functions depending on the domain and boundary conditions of the PDE's; e.g., Fourier basis on a periodic domain. The nonparametric nature can be understood as follows. If the diffusion processes in \eqref{SDE} are exactly the isotropic gradient flows \eqref{gradflow}, then $\cal L = \hat{\cal L}$ and, 
\BEA\index{Monte-Carlo approximation}
c_j(t+\tau) = \sum_l \langle \varphi_l,e^{\tau\mathcal{L}}\varphi_j\rangle_{\peq} c_l(t) = \sum_l e^{\lambda_l\tau} c_l(t) \approx \sum_{l=1}^M \sum_{i=1}^N e^{\lambda_l\tau} p_t(\theta_i)\varphi_l(\theta_i)\peq(\theta_i)^{-1},\label{coeff1}
\EEA
where $\lambda_l$ are eigenvalues of $\hat{\cal L}$ such that $\hat{\cal L}\varphi_l = \lambda_l\varphi_l$ and Monte-Carlo approximation (evaluated on data set $\theta_i$) is used to approximate the inner-product $\langle \cdot,\cdot \rangle$ defined with respect to $L^2(\mathcal{M})$. Here, $M$ denotes the number of eigenfunctions that are used in the numerical approximation and if $M$ is too small, then the Gibbs phenomena\index{Gibbs phenomena} will reduce the accuracy of the approximation as in the standard spectral method. With \eqref{forecast} and \eqref{coeff1}, the forecasting problem for gradient flows with isotropic diffusion can be solved without needing to know the expressions for $a, b,$ or $\hat{\cal L}$ and this is what we meant by nonparametric modeling. 

For general diffusion processes, we can approximate the semigroup solutions of the corresponding generator $\mathcal{L}$ with a shift operator\index{shift operator} defined as follows
\BEA
Sf(\theta_{i}) =  f(\theta_{i+1}), \label{shiftop}
\EEA
for any smooth function $f\in L^2(\mathcal{M},\peq)$. It was shown in \cite{bgh:15} that the stochastic operator $S$ is an unbiased estimator of $e^{\tau\cal L}$. Furthermore, the error from the stochastic nature of $S$ is minimized by representing $S$ in the diffusion basis coordinate $\varphi_j$ (eigenfunctions of $\hat{\cal L}$), see \cite{bgh:15} for details. As mentioned above, this is the second motivation for projecting the probabilistic forecasting problem in \eqref{FokkerPlanck} on these coordinate basis. In this general case (non-gradient drift anisotropic diffusions), $A_{jl}=\langle\varphi_l,e^{\tau\mathcal{L}}\varphi_j\rangle_{\peq}\approx \langle\varphi_l,S\varphi_j\rangle_{\peq}=\hat{A}_{jl}$ and the coefficients in \eqref{coeff1} become,
\BEA
c_j(t+\tau) &=& \sum_l A_{jl} c_l(t) \approx \sum_l \hat{A}_{jl} c_l(t),\nonumber\\
\hat{A}_{jl} &\approx& \frac{1}{N}\sum_{i=1}^N  \varphi_l(\theta_i)\varphi_l(\theta_{i+1}),
\EEA
using the definition of shift operator in \eqref{shiftop} and by Monte-Carlo averaging. For longer time, we can iterate $A$ to obtain $\vec{c}(t+n\tau)=A^n\vec{c}(t)$, where $c_j$ is the $j$-th component of $\vec{c}$. We should note that by the ergodicity assumption on the diffusion process in \eqref{SDE}, we ensure that the largest eigenvalue of $e^{\tau\mathcal{L}}$ is equal to 1 with constant eigenfunction, $\mathbbm{1}(\theta)$, and it can be shown that the largest eigenvalue of $A$ is also 1 corresponding to eigenvector $[\vec{e}_1]_j=\langle \mathbbm{1}, \varphi_j\rangle$, i.e., $A\vec{e}_1 = \vec{e}_1$. Here, $\vec{e}_1$ denotes a vector that is $1$ on the first component and zero otherwise. Therefore,  
\BEA
\lim_{n\rightarrow \infty} \vec{c}(t+n\tau) = \lim_{n\rightarrow \infty} A^n \vec{c}(t) = \vec{e}_1.
\EEA
 This means the forecast in \eqref{forecast} will converge to the equilibrium density, 
\BEA
\lim_{t\rightarrow \infty} p(\theta,t) =  \lim_{n\rightarrow\infty} \sum_j c_j (t+n\tau) \varphi_j(\theta)\peq(\theta) = \mathbbm{1}(\theta)\peq(\theta) = \peq(\theta).\label{consistent}
\EEA
Numerically, the largest eigenvalue of $\hat{A}$ can be greater than 1 due to finite samples in the Monte-Carlo integral. To overcome this issue, one can ensure the stability by dividing any eigenvalue with norm greater than 1 so it has norm equal to 1. Subsequently, the nonparametric forecast will produce a consistent equilibrium density as shown in \eqref{consistent}, by design.

\begin{figure*}
\centering
\includegraphics[width=0.7\textwidth]{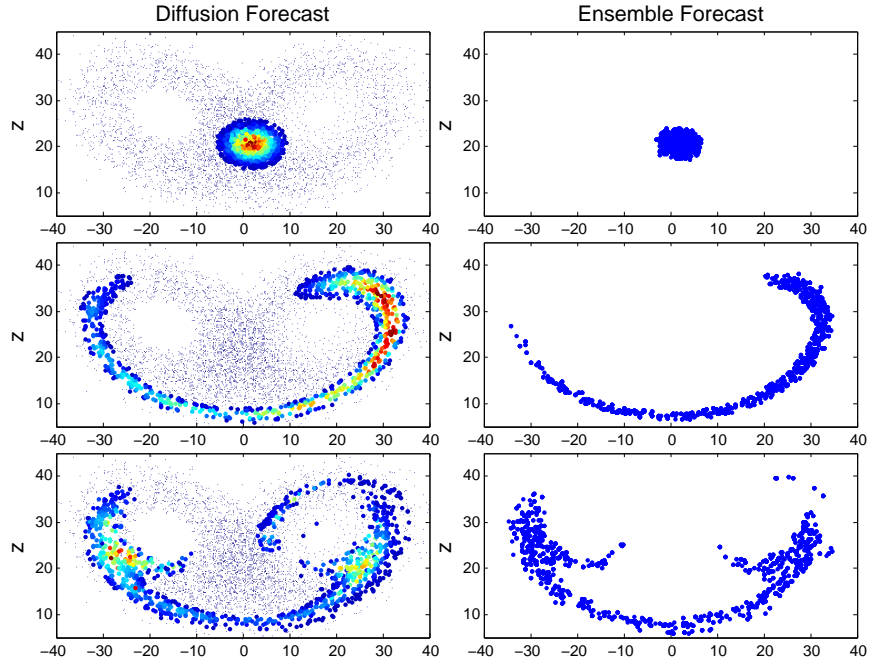}
\caption{Probability densities (as functions of $x+y$ and $z$) from the equation-free Diffusion Forecasting model (left column) and an ensemble forecasting (right column) at times $t=0$ (first row), $t=0.5$ (second row), and $t=2$ (third row). On the left column, the color spectrum ranging from red to blue is to denote high to low value of density.}
\label{L63fig2} 
\end{figure*}

\begin{example}[Example 6:]
In Figure~\ref{L63fig2}, we show snapshots of probabilistic density at various times, obtained from the equation-free, diffusion forecasting method, on the famous chaotic dynamical system, the three-dimensional Lorenz-63 model \cite{lorenz:63}.\index{Lorenz-63 model}\footnote{This example is taken from \cite{bgh:15}.} For comparison, we also show the Monte-Carlo approximation of the evolution of the density (or ensemble forecasting), assuming that the full Lorenz-63 model is known. In this experiment, the same Gaussian initial conditions are prescribed (as shown in the first row in Figure~\ref{L63fig2}). In each panel of this figure, we show the density as functions of $x+y$ and $z$ (corresponding to the three components of the Lorenz model). In the left column, we also show the data set that are used for training the diffusion model (smaller black dots). Notice that even at a long time $t=2$ (which is longer than the doubling time of this model, 0.78), the densities obtained from both forecasting methods are still in a good agreement.
\end{example}

\subsection{Semiparametric framework to mitigate model error}
In the nonparametric framework above, the diffusion forecasting model \emph{interpolates} from the training data, meaning the required data must fill in the manifold. This implies that the required data grows exponentially as a function of the dimension of the manifold $\mathcal{M}$ and this is the practical limitation of such approach. For high-dimensional systems, however, we usually have some physical knowledge, but these high-dimensional parametric models are subject to model error as discussed before. The idea of semiparametric modeling is to use the nonparametric model to compensate for the low-dimensional model error in the high-dimensional parametric model. 
 
In the current semiparametric framework \cite{bh:15qjrms}, we assume that the underlying truth solves,\index{semiparametric framework}
\BEA
\dot x &=& f(x,\theta),\nonumber\\
d\theta &=& a(\theta) \, dt + b(\theta)\, dW_t, \nonumber
\EEA
where the parametric model $f$ is known but neither the parameter $\theta$ nor its dynamics, $a, b$, are known. Here, we are assuming that model error is attributed to imperfect specification of dynamically evolving parameters $\theta$ with unknown dynamics. Presently, there are more implicit assumptions for such framework to work, including $\theta\in\mathcal{M}$ to be low dimensional and independent of $x$. While the high dimensionality issue will still be the fundamental practical issue for this framework, the second issue, constructing diffusion forecasting models for conditionally distributed data $\theta_i\sim p(\theta|x)$ is an important open problem that we plan to address in near future. With all these assumptions, let us demonstrate the semiparametric framework for mitigating model error.

The first step is to extract a time series of $\theta$ from noisy observations, 
\BEA
v_m = h(x_m) + \epsilon_m, \quad \epsilon_m\sim\mathcal{N}(0,R).\label{finalobs}
\EEA
Obviously, when function $h$ also depends on $\theta$, this problem becomes simpler assuming that the theoretical observability condition is satisfied; but in most applications, $\theta$ is hidden and we still assume that  the theoretical observability condition is satisfied as in any standard inverse problems. In \cite{bh:15qjrms}, we demonstrate that the time series for $\theta$ can be extracted by implementing the adaptive covariance estimation method discussed in Section~\ref{sec22}, treating $\theta$ as Gaussian white noise processes. Using the extracted training data set, we build a nonparametric model for $p(\theta,t)$ with the strategy discussed in Section~\ref{sec51}.  Subsequently, we combine the parametric and nonparametric models by sampling $\theta^k(t) \sim p(\theta,t)$ from the nonparametric model to be used with the ensemble forecast $(x^k,\theta^k)$, where subscript $k$ denotes the $k$-th ensemble member. Again, see \cite{bh:15qjrms} for the implementation detail. We should note that 
this framework was designed to maintain as much of the current parametric ensemble forecasting and filtering framework (as used in the numerical weather prediction) as possible.

\begin{example}[Example 7:] 
Here, we demonstrate the application of the semiparametric framework on the following system,
\BEA
\frac{dx_j}{dt} &=& \theta x_{j+1}x_{j-1}-x_{j-2}x_{j-1} - x_j + 8,\nonumber\\
\theta &=& \frac{x}{40}+1 \nonumber \\
\dot x &=&10(y - x),  \label{semiparametric}\\
\dot y &=& 28 x - y - xz,\nonumber \\
\dot z &=& xy - \frac{8}{3}z, \nonumber
\EEA
where the parameter $\theta$ is a rescaling of $x$ that is dynamically evolving in accordance to the Lorenz-63 model \cite{lorenz:63}. The rescaling is to confine $\theta\in[0.5,1.5]$ to avoid numerical instability since the quadratic terms do not conserve the energy-like quantity, $E=\sum_j x_j^2$, when $\theta\neq 1$. So, bad estimates of $\theta$ beyond this interval can produce unstable forecasts. 

For comparison, we include forecasting results from the perfect model (which assumes knowing the full system in \eqref{semiparametric}) and from Lorenz-96 model in \eqref{L96} (this is equivalent to assuming $\theta=1$ in \eqref{semiparametric} and ignoring the Lorenz-63 model). For a more complete comparison with other approaches such as the Heterogeneous-Multiscale-Methods \cite{eelrv:07}\index{HMM (Heterogeneous Multiscale Methods)} or persistence model, see \cite{bh:15qjrms}. In each of these experiments, the ensemble forecasts are performed with 86 ensemble members (doubling the total state variables of the full model). Notice that the perfect model produces the best short term prediction, but it also seems to produce a biased forecast (based on larger RMS error above,  the climatological error, in the intermediate time beyond 12 days). We suspect that this bias is due to the sampling error introduced by finite ensemble size. On the other hand, the RMS error of the semiparametric forecast grows slightly more quickly than that of the perfect model initially, but the forecast is unbiased in the intermediate time, approaching the climatology without exceeding the climatological error. After 8 model days the semiparametric forecast produces a better forecast compare to the perfect model. We suspect that this is  because the samples of the nonparametric forecast are independent to the forecast density. Finally, the standard Lorenz-96 model produces the worse forecast.
\end{example}
\begin{figure*}
\centering
\includegraphics[width=0.5\textwidth]{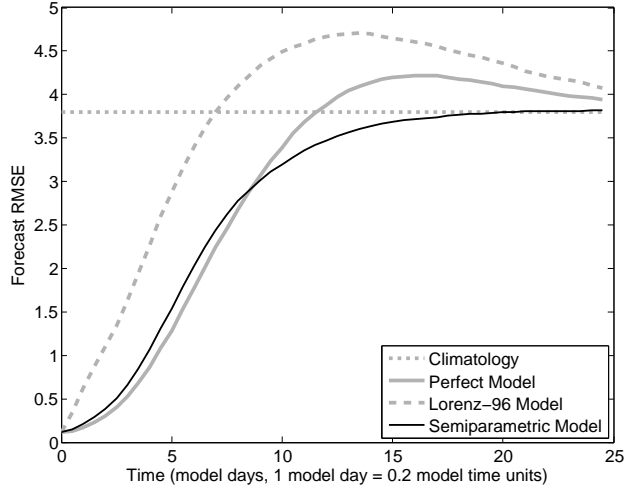}
\caption{Comparison of forecasting errors as functions of time.}
\label{L63} 
\end{figure*} 

While this result is encouraging as a first conceptual proof that it is possible to use the nonparametric modeling approach as a model error estimator, there are still many open questions related to various assumptions that are limiting the application of this framework on more complex problems.

\section{Summary}
In this chapter, we discussed one major challenge in data assimilation: model error. Various perspectives of model error were offered. First, the traditional point of view, which is based on a prior distribution formulation was overviewed and compared with a more recent, posterior distribution, formulation. Simple examples were used to elucidate the robustness of the posterior distribution formulation. Second, various methods to mitigate model error were discussed. We classified these methods into two categories: the \emph{statistical methods} for those who directly estimate the low-order model error statistics; and the \emph{stochastic parameterizations} for those who implicitly estimate all statistics by imposing stochastic models beyond the traditional unbiased white noise Gaussian processes. We hope that this discussion also clarifies a common misconception in the data assimilation community of associating model error to only estimating the model error covariance, $Q$. Indeed, the posterior distribution formulation shows that even in simple contexts, the optimal model error estimator involves parametric terms that depend on the estimates from an imperfect model. Third, for model error due to unresolved scales, connection to related subjects under different names in applied mathematics, such as the Mori-Zwanzig formalism and the averaging method, were discussed with the hope that the existing methods can be more accessible and eventually be used appropriately.
Fourth, we provide a theoretical foundation to support the use of stochastic parameterization to mitigate model error in data assimilation and point out the fundamental issues in lifting this approach for general problems. Namely, the difficulties in choosing the appropriate (stable and consistent) models and in designing efficient and accurate schemes to estimate the parameters in the parametric models. Fifth, we show an alternative strategy for mitigating model error with a nonparametric approach, using stable and consistent data-driven models constructed with diffusion maps algorithms. While the idea works under various assumptions, we hope that this result motivates the development in this direction to handle more complex problems. 

While covariance inflation with an empirically chosen $Q$ matrix has been the most popular approach in mitigating model error in data assimilation since it is the most practical numerically, we hope that the review in this chapter can provide more compelling reasons for alternative approaches such as the stochastic parameterization and nonparametric modeling to be seriously considered. A significant challenge remains, that is, to apply all these theoretically profound techniques on realistic, large-scale applications in which most assumptions are violated or unverifiable. To make progress, more interdisciplinary collaborative effort is crucial. In particular, we advocate for synergistic collaborative efforts between physicists for their physical intuitions, mathematicians for theoretical justifications, and engineers for efficient implementations. This chapter is written with the hope that it provides a connection between theoreticians and practitioners for such collaboration. In particular, many discussions regarding to advantages and limitations on various methods in this chapter can be useful to motivate more interdisciplinary research in this field.

\section*{Acknowledgment} The author research is partially supported by the Office of Naval Research Grants N00014-11-1-0310, N00014-13-1-0797, MURI N00014-12-1-0912 and the National Science Foundation DMS-1317919. The author thanks J.L.~Anderson, F. Zhang, A.J. Majda, T. Berry, and anonymous readers for their valuable inputs and criticisms that help improving this note. 


\printindex

\end{document}